\documentclass[a4paper,11pt]{amsart}

\newfont{\cyr}{wncyr10 scaled 1100}

\usepackage[left=2.7cm,right=2.7cm,top=3.5cm,bottom=3cm]{geometry}
\usepackage{amsthm,amssymb,amsmath,amsfonts,mathrsfs,amscd,yfonts}
\usepackage{verbatim}
\usepackage{turnstile}
\usepackage[latin1]{inputenc}
\usepackage[all]{xy}
\usepackage{latexsym}
\usepackage{stmaryrd}
\usepackage{dsfont}
\usepackage{longtable}
\usepackage{multirow}
\usepackage[usenames,dvipsnames]{color}
\usepackage{listings}
\usepackage{mathtools}
\usepackage{tikz}
\usepackage{tikz-cd}
\usetikzlibrary{matrix,arrows,decorations.pathmorphing}
\usepackage{bm}
\usepackage[shortlabels]{enumitem}
\usepackage[pagebackref]{hyperref}

\theoremstyle{plain}

\newtheorem{theorem}{Theorem}[section]
\newtheorem{corollary}[theorem]{Corollary}
\newtheorem{lemma}[theorem]{Lemma}
\newtheorem{proposition}[theorem]{Proposition}
\newtheorem{propo}[theorem]{Proposition}

\newtheorem{conj}[theorem]{Conjecture}

\theoremstyle{definition}
\newtheorem{definition}[theorem]{Definition}
\newtheorem{defi}[theorem]{Definition}

\newtheorem{examplewr}[theorem]{Example}
\newtheorem{ass}[theorem]{Assumption}

\theoremstyle{remark}
\newtheorem{obswr}[theorem]{Observation}
\newtheorem{remarkwr}[theorem]{Remark}

\newenvironment{remark}{\begin{remarkwr}\begin{upshape}}{\end{upshape}\end{remarkwr}}

\newcommand{\bb}{\mathbb}
\newcommand{\frk}{\mathfrak}
\newcommand{\cl}{\mathcal}

\DeclareMathOperator{\BK}{BK}
\DeclareMathOperator{\BF}{BF}

\DeclareMathOperator{\cris}{cris}

\DeclareMathOperator{\et}{et}

\DeclareMathOperator{\cyc}{cyc}

\DeclareMathOperator{\rec}{rec}

\DeclareMathOperator{\Ind}{Ind}

\DeclareMathOperator{\lcm}{lcm}
\DeclareMathOperator{\bal}{bal}
\DeclareMathOperator{\Eis}{Eis}

\DeclareMathOperator{\Log}{Log}

\DeclareMathOperator{\Frac}{Frac}
\DeclareMathOperator{\ac}{ac}
\DeclareMathOperator{\ab}{ab}

\newcommand{\Lp}{{\mathscr{L}_p}}

\newcommand{\Q}{\mathbb{Q}}
\newcommand{\Z}{\mathbb{Z}}

\newcommand{\Sel}{\mathrm{Sel}}

\newcommand{\Gal}{\mathrm{Gal\,}}
\newcommand{\GL}{\mathrm{GL}}

\newcommand{\Fr}{\mathrm{Fr}}

\newcommand{\ord}{{\mathrm{ord}}}
\newfont{\gotip}{eufb10 at 12pt}

\newcommand{\hE}{{\mathbf{E}}}
\newcommand{\hf}{{\mathbf{f}}}
\newcommand{\hg}{{\mathbf{g}}}
\newcommand{\hh}{{\mathbf h}}

\newcommand{\res}{\mathrm{res}}

\numberwithin{equation}{section}

\include{thebibliography}

\begin{document}

\title[Anticyclotomic diagonal cycles and Beilinson--Flach]{Anticyclotomic diagonal classes and Beilinson--Flach elements}

\author{Ra\'ul Alonso, Lois Omil-Pazos and \'Oscar Rivero}

\begin{abstract}
We present a comparison between the anticyclotomic Euler system of diagonal cycles associated with the convolution of two modular forms and the cyclotomic Beilinson--Flach Euler system. This extends the seminal work of Bertolini, Darmon, and Venerucci, who established a link between (anticyclotomic) Heegner points and the Beilinson--Kato system. Our approach hinges on a detailed analysis of $p$-adic $L$-functions and Perrin-Riou maps and exploits the Eisenstein degeneration of diagonal cycles along Hida families, working with a CM family which specializes to an irregular Eisenstein series in weight one. We use these results to derive some arithmetic applications.
\end{abstract}

\date{\today}

\address{R. A.: School of Mathematics and Statistics, University College Dublin, Belfield, Dublin 4, Ireland}
\email{raul.alonsorodriguez@ucd.ie}

\address{L. O.-P.: CITMAga, Santiago de Compostela, Spain}
\email{lois.omil.pazos@usc.es}

\address{O. R.: CITMAga and Departamento de Matem\'aticas, Universidade de Santiago de Compostela, Santiago de Compostela, Spain}
\email{oscar.rivero@usc.es}


\maketitle

\setcounter{tocdepth}{1}
\tableofcontents

\section{Introduction}

The recent work of Bertolini, Darmon, and Venerucci \cite{BDV} establishes a conjecture relating Heegner cycles to Beilinson--Kato elements, linking both objects to $p$-adic families of Beilinson--Flach elements in the higher Chow groups of products of two modular curves. The comparison between Heegner points and Beilinson--Kato elements had previously been investigated by B\"uy\"ukkboduk \cite{buyukboduk16} and Venerucci \cite{Ven}, within the framework of the exceptional zero conjectures. In this work, we extend the results of \cite{BDV} by relating Beilinson--Flach classes to the diagonal cycles studied in \cite{BSV} and \cite{DR3}, which in turn can be interpreted in terms of the anticyclotomic classes constructed in \cite{ACR}, \cite{ACR2}, and \cite{castella-do}.

This work may also be interpreted in the context of the recent developments of Loeffler and the third author \cite{LR1}, where they study the Eisenstein degeneration of Euler systems, i.e. Euler systems obtained from critical Eisenstein series. Our setting, however, differs in a significant way, since it involves as input a (non-cohomological) weight-one Eisenstein series. We will return to this point throughout the text.

Following \cite{BeDiPo}, the central feature in the setting of irregular weight-one modular forms attached to Dirichlet characters $(\chi_1,\chi_2)$ with $\chi_1(p)=\chi_2(p)$ is the existence of three (ordinary) families passing through it: the families of Eisenstein series $\hE(\chi_1,\chi_2)$ and $\hE(\chi_2,\chi_1)$, as well as a third Hida family that is generically cuspidal. The key idea, both in this note and in earlier related works, is the introduction of an auxiliary modular form that can be suitably deformed along a Hida or Coleman family to produce the required Galois representation. In this sense, the present work may be viewed as an instance of CM degeneration, while the work of Loeffler and the third author provides an example of Eisenstein degeneration.

\subsection{The set-up}

Let $p$ be an odd prime and fix once and for all an embedding $\iota_p:\overline{\bb{Q}}\xrightarrow{} \overline{\bb{Q}}_p$ and an embedding $\iota: \overline{\bb{Q}}\hookrightarrow \bb{C}$. The latter determines a choice of complex conjugation in $G_{\bb{Q}}$ which will be denoted by $c$.

Let $K/\Q$ be an imaginary quadratic field in which $p$ splits and let $\varepsilon_K$ be the quadratic character attached to it. Let $\hf$ be the CM Hida family introduced in \S\ref{subsec:wt1}, which specializes in weight one to the irregular weight-one Eisenstein series $f = \Eis_1(\varepsilon_K)$. Let $(\hg,\hh)$ be a pair of Hida families of coprime tame levels $(N_g,N_h)$ and characters $(\chi_\hg,\chi_\hh)$ such that $\chi_\hg \chi_\hh =\varepsilon_K \omega^{2r}$ for some $r\in \bb{Z}$. We also assume that $\hg$ and $\hh$ are residually irreducible and $p$-distinguished. Let $g=\hg_{y_0}^\circ$ and $h=\hh_{z_0}^\circ$ be good crystalline specializations of $\hg$ and $\hh$ of weights $l_0\geq 2$ and $m_0\geq 1$, respectively. We assume that $p\nmid \mathrm{cond}(h)$. Set $c_0=(l_0+m_0-1)/2$.


In \S\ref{subsec:families} we introduce our notations regarding Galois representations attached to families of modular forms. In particular, attached to $\hg$ (and similarly for $\hf$ and $\hh$) there is a locally-free rank-two module $\mathbb V_{\hg}$, defined over a finite flat extension of $\Lambda = \mathbb Z_p[[1+p\mathbb Z_p]]$ and equipped with a continuous action of the absolute Galois group $G_{\mathbb Q}$.

\subsection{Euler systems and $p$-adic $L$-functions}

We now present the two cohomology class that can naturally be attached to the triple $(f,\hg,\hh)$. Firstly, the diagonal cycle class, which may be understood as an anticyclotomic Euler system as discussed in \cite{ACR}; secondly, the cyclotomic system of Beilinson--Flach, as developed e.g. in \cite{KLZ}.

\begin{enumerate}
    \item[(A)] The work of Bertolini--Seveso--Venerucci \cite{BSV} and Darmon--Rotger \cite{DR3} provides us with a diagonal cycle class \[ \kappa(\hf,\hg,\hh) \in H^1_{\bal}(\mathbb Q,  \bb{V}_{\hf} \hat\otimes \mathbb V_{\hg} \hat \otimes \mathbb V_{\hh}(2-{\bf t})), \] where $2{\bf t} = {\bf k} + {\bf l} + {\bf m}$. This class is constructed via the $p$-adic variation of diagonal cycles, so it may be understood as a geometric object. By specilizing the first variable, it yields a class $\kappa(f,\hg,\hh)\in H^1_{\bal}(\bb{Q},V_f\otimes\bb{V}_{\hg}\hat{\otimes}\bb{V}_{\hh}(2-\mathbf{t}_1))$, where $2\mathbf{t}_1=1+\mathbf{l}+\mathbf{m}$.

    \item[(B)] The work of Kings--Loeffler--Zerbes \cite{KLZ} provides us with a (cyclotomic) Beilinson--Flach class
    \[
    \kappa_{\hg,\hh}\in H^1_{\bal}(\mathbb Q, \mathbb V_{\hg} \hat \otimes \mathbb V_{\hh} (2-\mathbf{t}_1))
    \]
    attached to the pair $(\hg,\hh)$, with the conventions that we later recall. Similarly, we may consider a Beilinson--Flach class $\kappa_{\hg,\hh\otimes\varepsilon_K}$ attached to the families $(\hg, \hh \otimes \varepsilon_K)$, where $\hh\otimes \varepsilon_K$ denotes the twist of $\hh$ by the quadratic character $\varepsilon_K$.
\end{enumerate}

Note, however, that the construction of the Beilinson--Flach classes proceeds in an ostensibly different way, since it involves working with modular units, which are absent in the theory of diagonal cycles. Roughly speaking, the construction of diagonal cycles proceeds in purely geometric terms, so the construction is amenable to be generalized to other settings where there are no modular units (e.g. Shimura curves) and Beilinson--Flach classes are not available.

The main result of this note is a result connecting both classes. However, to ensure that they live in the same space, one first defines a Beilinson--Flach class \[ \BF(f,\hg,\hh) \in H^1(\mathbb Q, V_f \otimes \mathbb V_{\hg} \hat{\otimes} \mathbb V_{\hh}(2-{\mathbf t}_1)), \] obtained as a suitable weighted combination of $\kappa_{\hg,\hh}$ and $\kappa_{\hg,\hh \otimes \varepsilon_K}$. The intuition for that comes from the observation that $V_f \otimes \mathbb V_{\hg} \hat \otimes \mathbb V_{\hh}(2-{\mathbf t}_1)$ decomposes as the direct sum of the Galois representations $\mathbb V_{\hg} \hat \otimes \mathbb V_{\hh}(2-{\mathbf t}_1)$ and $\mathbb V_{\hg} \hat \otimes \mathbb V_{\hh}(2-{\mathbf t}_1)(\varepsilon_K)$, and hence one can naturally construct a cohomology class by gluing the Euler systems for each of the two pieces.

The key point for developing a comparison of Euler systems comes through the theory of $p$-adic $L$-functions. More precisely, the arithmetic of the Beilinson--Flach system is dictated by the Hida--Rankin $p$-adic $L$-functions, which interpolate special values of the complex $L$-function $L(g \otimes h,s)$. Similarly, diagonal cycles are related, via explicit reciprocity laws, with the triple product $p$-adic $L$-function constructed by Hsieh \cite{Hs}, and later extended by Andreatta--Iovita \cite{AI} to the non-ordinary case.

The main result of this note may be interpreted as a comparison between a cyclotomic and an anticyclotomic cohomology class attached to the representation $V_f \otimes \mathbb V_{\hg} \hat \otimes \mathbb V_{\hh}(2-{\mathbf t}_1)$. As a piece of notation, we put $\bb{V}_{f\hg\hh}^\dagger=V_f\otimes\bb{V}_{\hg}\hat{\otimes} \bb{V}_\hh (2-\mathbf{t}_1)$ and $\bb{V}_{\hg\hh}^\dagger=\bb{V}_{\hg}\hat{\otimes}\bb{V}_\hh(2-\mathbf{t}_1)$. Before stating the theorem, we introduce the following objects; their precise definitions are recalled in the main body of the article.

\begin{itemize}
    \item[-] The Selmer group $H^1_{\cl{G}\cup +}(\bb{Q}, \bb{V}_{f\hg\hh}^\dagger)$, corresponding to a certain local condition at $p$. This is defined in  \S\ref{subsec:Selmergroups}, and it is a module over $\Lambda_{\hg\hh}=\Lambda_{\hg}\hat{\otimes}\Lambda_{\hh}$, where $\Lambda_{\hg}$ (resp. $\Lambda_{\hh}$) is the Iwasawa algebra over which the Hida family $\hg$ (resp. $\hh$) is defined.
    \item[-] The Selmer group $H^1_{\bal}(\bb{Q}, \bb{V}_{f\hg\hh}^\dagger)$, corresponding to the balanced local condition at $p$.
    \item[-] The Atkin--Lehner pseudo-eigenvalue $\lambda_{N_g}(\hg)$.
    \item[-] The $p$-adic period $\Omega_{f,\gamma}$, depending on the choice of the isomorphism $\gamma$ of \eqref{eq:iso}.
    \item[-] The triple product $p$-adic $L$-function $\Lp^g(\hf,\hg,\hh)$, interpolating central values of the corresponding triple product complex $L$-functions along the $\hg$-unbalanced region. In this case, $\Lp^g(f,\hg,\hh)$ corresponds to specializing $\hf$ at $f$.
\end{itemize}

In the course of developing the theory and proving the results of this work, we impose a non-vanishing assumption on $\Lp^g(f,\hg,\hh)$; see Assumption \ref{ass:nonvanishing}.

\begin{theorem}
Assume that $H^1_{\cl{G}\cup +}(\bb{Q}, \bb{V}_{f\hg\hh}^\dagger)$ is a torsion-free $\Lambda_{\hg\hh}$-module of rank $1$ and that the $p$-adic $L$-function $\Lp^g(f,\hg,\hh)$ is not identically zero. Then $\BF(f,\hg,\hh)$ belongs to $H^1_{\bal}(\bb{Q}, \bb{V}_{f\hg\hh}^\dagger)$ and \[ \lambda_{N_g}(\hg)\cdot\BF(f,\hg,\hh)=\Omega_{f,\gamma}\cdot \Lp^g(f,\hg,\hh) \cdot \kappa(f,\hg,\hh). \]
\end{theorem}

Furthermore, the result may be also understood in the framework of the comparison between different instances of Euler systems, developed e.g. in \cite{LR1}, but where one considers instead CM families passing through a weight-one Eisenstein series, while in \emph{loc.\,cit.} the key input was the use of families passing through the critical $p$-stabilization of an Eisenstein series of weight at least two.

\begin{remark}
Contrary to the situation in \cite{BDV}, where Beilinson--Kato classes and Beilinson--Flach elements were compared in Iwasawa cohomology, we are not working at that level here. In our case, such a comparison is not possible, since the main result should be viewed as relating a cyclotomic Euler system to an anticyclotomic one. Therefore, we restrict our analysis to the bottom layers.
\end{remark}

\begin{remark}
The slogan of this result may be summarized under the sentence {\it Heegner points are to Kato classes what anticyclotomic diagonal cycles are to Beilinson--Flach classes}. Unfortunately, this setting remains much more mysterious in many instances, like the lack of a proof of the Iwasawa main conjecture (only one divisibility is known).
\end{remark}

The proof proceeds in three main steps.

\begin{enumerate}
\item[(a)] Compare the triple product $p$-adic $L$-function $\Lp^g(f,\hg,\hh)$ with the Hida--Rankin $p$-adic $L$-functions attached to the pairs $(\hg,\hh)$ and $(\hg, \hh \otimes \varepsilon_K)$.
\item[(b)] Define the weighted Beilinson--Flach class and show that it satisfies the appropriate local condition. This requires a careful analysis of the structure of the corresponding Selmer groups.
\item[(c)] Relate the two cohomology classes via the explicit reciprocity laws.
\end{enumerate}

As a consequence of this result in families, we also obtain a comparison upon specializing $\hg$ and $\hh$; see \S\ref{subsec:nonexceptional} for details. Another corollary concerns a factorization of a certain big logarithm, which can be applied to the Beilinson--Flach class. To state it, we introduce the following objects:
\begin{itemize}
    \item[-] The triple product $p$-adic $L$-function of Hsieh \cite{Hs} and \cite{AI}, denoted $\Lp^f(\hf,\hg,\hh)$, attached to the triple $(\hf,\hg,\hh)$, where $\hf$ is the CM family considered above. Its specialization at weight one in the first variable is denoted by $\Lp^f(f,\hg,\hh)$.
    \item[-] The big-logarithm map $\Log_{\omega_{\hg} \otimes \omega_{\hh}}$, introduced in Def. \ref{def:big-log}.
\end{itemize}

\begin{corollary}
Assume that $H^1_{\cl{G}\cup +}(\bb{Q}, \bb{V}_{f\hg\hh}^\dagger)$ is a torsion-free $\Lambda_{\hg\hh}$-module of rank $1$ and that the $p$-adic $L$-function $\Lp^g(f,\hg,\hh)$ is not identically zero. Then we have the following factorization of the image of the Beilinson--Flach class under the Perrin-Riou map $\Log_{\bm{\omega}_{\hg} \otimes \bm{\omega}_{\hh}}$: \[ \Omega_{f,\gamma} \cdot \Lp^g(f,\hg,\hh) \cdot \Lp^f(f, \hg, \hh) =  \lambda_{N_g}(\hg) \cdot \Log_{\omega_{\hg} \otimes \omega_{\hh}}(\BF(f, \hg, \hh)). \] 
\end{corollary}

\subsection{Exceptional zeros}

There is a case which is especially interesting from the point of view of arithmetic applications, which corresponds to the following assumption.

\begin{ass}
With the previous notations, it holds that \[ \alpha_g \beta_h = p^{\frac{l_0+m_0-3}{2}}. \] 
\end{ass}

Since $g$ and $h$ are $p$-ordinary and $p\nmid \mathrm{cond}(h)$, it follows from this assumption together with the Ramanujan--Petersson conjecture that $(l_0,m_0)=(2,1)$ and that $g$ has conductor $N_gp$. Note that this situation includes in particular the case where $g$ is the modular form attached to an elliptic curve over $\bb{Q}$ with multiplicative reduction at $p$.

Hence, we encounter an arithmetic situation that is particularly intriguing, as {\em two} of the Euler factors associated with the $p$-adic $L$-function vanish. In this setting, improved $p$-adic $L$-functions are available, both in the Hida--Rankin case and for the triple product of modular forms. However, the vanishing of these two Euler factors occurs for different reasons: one is introduced by the big logarithm map, while the other arises from the interpolation in families of Beilinson--Flach classes and diagonal cycles. To complete the picture, we therefore need to (a) construct improved big logarithm maps, and (b) construct improved cohomology classes. At this point, we rely on a standard conjecture in the theory of exceptional zeros to guarantee the existence of the improved Beilinson--Flach class in this setting, while the improved diagonal cycle class has already been constructed in \cite{BSV}.

As a piece of notation for the next proposition, we use the following terminology.
\begin{itemize}
    \item[-] The line $\mathcal C$ corresponds to the points in the pair of families with weights of the form $(l,l-1)$; see \S\ref{sec:expectional} for the precise definition.
    \item[-] The triple-product improved $p$-adic $L$-function, introduced in Theorem \ref{thm:improvedtripleproductpadicLfunction}, is denoted by $\widehat{\Lp}^g(f,\hg,\hh)$.
    \item[-] The improved diagonal cycle of \cite{BSV} is denoted by $\widehat{\kappa}(f,\hg,\hh)$.
    \item[-] The (still conjectural) improved Beilinson--Flach class is denoted by $\widehat{\mathrm{BF}}(f, \hg, \hh)$.
\end{itemize}

\begin{theorem}
Assume that $H^1(\bb{Q},\bb{V}_{f\hg\hh}^\dagger\vert_{\cl{C}})$ is a torsion-free $\cl{O}_{\cl{C}}$-module, $H^1_{\cl{G}\cup +}(\bb{Q},\bb{V}_{f\hg\hh}^\dagger)$ is a torsion-free $\Lambda_{\hg\hh}$-module of rank $1$ and $\Lp^g(f,\hg,\hh)$ is not identically zero. Then $\widehat{\mathrm{BF}}(f, \hg, \hh)$ belongs to the Selmer group $H^1_{\bal}(\bb{Q},\bb{V}_{f\hg\hh}^\dagger\vert_{\cl{C}})$ and
\[
\Omega_{f,\gamma} \cdot \widehat{\Lp}^g(f,\hg,\hh) \cdot \widehat{\kappa}(f,\hg,\hh) = \lambda_{N_g}(\hg)\cdot\widehat{\mathrm{BF}}(f, \hg, \hh).
\]
\end{theorem}

The proof follows the same ideas as in the non-exceptional case, but requires certain modifications to account for the improved setting. In particular, one needs to work with improved big logarithms, since in this situation the Euler factor in the numerator of the Perrin-Riou map vanishes at the point corresponding to the pair $(g,h)$.

\subsection{Related works}

This project, together with \cite{BDV}, provides an example of how two distinct types of Euler systems can be related through a {\it CM degeneration technique}. We now highlight its connections with other analogous phenomena:

\begin{itemize}
    \item[(1)] {\bf Critical Eisenstein series.} In \cite{LR1} and \cite{PR25}, the authors study Euler systems in families where one of the modular forms passes through a point corresponding to the critical-slope $p$-stabilization of an Eisenstein series. In these cases, the specialization of the associated Galois representation admits a projection onto a one-dimensional quotient, allowing for a direct comparison with a {\it smaller} Euler system.

    \item[(2)] {\bf Irregular Euler systems.} The degeneration technique used in this work, inspired by \cite{BSV}, considers a CM family intersecting an irregular weight-one Eisenstein series. This can be viewed as a {\it CM degeneration}, in which the Euler system in families decomposes into a sum of two distinct Euler systems, each associated with a different component of the CM representation. In both situations, the key idea is the same: to exploit a cuspidal Hida or Coleman family which, at a certain point, specializes to an Eisenstein series.
    
    \item[(3)] {\bf (Eisenstein) congruences among modular forms.} Instead of working inside a family, one may also consider an Euler system attached to a cuspidal modular form that satisfies a congruence relation with an Eisenstein series. In this setting, similar connections to those described in (1) can be established.  

\end{itemize}

\subsection{Acknowledgements}
We are deeply grateful to Francesc Castella for stimulating discussions related to the subject of this work. We also thank K\^az\i m B\"uy\"ukboduk for many insightful comments and suggestions. The third-named author further thanks David Loeffler and Victor Rotger for numerous conversations that inspired several of the ideas developed in this note.

During the preparation of this article, the authors were supported by PID2023-148398NA-I00, funded by MCIU/AEI/10.13039/501100011033/FEDER, UE. The first author was also supported by Taighde \'Eireann -- Research Ireland under Grant number IRCLA/2023/849 (HighCritical).

\section{Preliminaries}

Along this section we recall the main properties about the Galois representations attached to modular forms and families which are needed along this work. Further, we state the assumptions which are needed to guarantee the \'etaleness of the Coleman--Mazur--Buzzard eigencurve around a classical point, following the seminal works of \cite{bellaiche-dimitrov} and \cite{BeDiPo}. This has a special relevance in our study, since we are going to consider a family passing through an irregular weight-one Eisenstein series.

Let $p>2$ be a prime and fix once and for all an embedding $\iota_p:\overline{\bb{Q}}\xrightarrow{} \overline{\bb{Q}}_p$ and an embedding $\iota: \overline{\bb{Q}}\hookrightarrow \bb{C}$. The latter determines a choice of complex conjugation in $G_{\bb{Q}}$ which will be denoted by $c$.

\subsection{Deligne representations}

Let $\xi=\sum_{n=1}^\infty a_n(\xi)q^n\in S_k(N_\xi,\chi_\xi)$ be a normalized newform of weight $k\geq 2$, level $N_\xi$, and nebentypus $\chi_\xi$. 
Let $L$ be a finite extension of $\bb{Q}_p$ containing the Fourier coefficients of $\xi$ (under the embedding $\iota_p$ fixed above) and let $\cl{O}$ be the ring of integers of $L$. By work of Eichler--Shimura and Deligne, there is a two-dimensional representation
$$
\rho_\xi\,:\,G_\mathbb Q\longrightarrow \GL_L(V_\xi)\simeq\GL_2(L)
$$
unramified outside $pN_\xi$ and characterized by the property
$$
{\rm trace}\,\rho_\xi(\Fr_{q})=a_{q}(\xi)
$$
for all primes $q\nmid pN_\xi$, where $\Fr_{q}$ denotes an arithmetic Frobenius element at $q$.
Let $Y_1(N_\xi)$ be the open modular curve over $\mathbb Q$ parameterizing pairs $(A,P)$ consisting of an elliptic curve $A$ and a point $P\in A$ of order $N_\xi$. Let $\mathscr{L}_{k-2}$ be the $p$-adic sheaf over $Y_1(N_\xi)$ introduced in \cite[\S2.3]{BSV}. We shall work with the geometric realization of $V_\xi$ arising as the maximal quotient of
\[
H^1_{\et}(Y_{1}(N_\xi)_{\overline{\mathbb Q}},\mathscr{L}_{k-2}(1))\otimes_{\Z_p}L
\]
on which the dual Hecke operators $T_q'$ and $\langle d\rangle'$ act as multiplication by $a_q(\xi)$ and $\chi_\xi(d)$ for all primes $q\nmid N_\xi$ and all $d\in(\Z/N_\xi\Z)^\times$. 



If $\xi$ is $p$-ordinary, then we have a filtration
\[
0\rightarrow V_\xi^+\rightarrow V_\xi \rightarrow V_\xi^-\rightarrow 0
\]
of $G_{\bb{Q}_p}$-representations, where $V_\xi^{\pm}$ are $1$-dimensional and $V_\xi^-$ is unramified with $\Fr_p$ acting as multiplication by the unit root of the $p$-th Hecke polynomial of $\xi$.


\subsection{Hida families and big Galois representations}\label{subsec:families}

In this subsection, we recall the fundamental concepts and main notations related to Hida families that will appear in this work.

Let $\Lambda=\bb{Z}_p[[1+p\bb{Z}_p]]$. Given an integer $r$ and a finite-order character $\epsilon:1+p\bb{Z}_p\rightarrow \overline{\bb{Q}}_p^\times$, we define the character $\nu_{r,\epsilon}:1+p\bb{Z}_p\rightarrow \overline{\bb{Q}}_p^\times$ by $\nu_{r,\epsilon}(z)=z^r\epsilon(z)$. This character extends to a ring homomorphism $P_{r,\epsilon}:\Lambda \rightarrow \overline{\bb{Q}}_p$, and we will use the same notation $P_{r,\epsilon}$ to denote the corresponding point in $\mathrm{Spec}(\Lambda)(\overline{\bb{Q}}_p)$. Points in $\mathrm{Spec}(\Lambda)(\overline{\bb{Q}}_p)$ of this form will be called \textit{arithmetic points}, and the set of such points will be denoted by $\cl{W}$. Given a finite flat extension $\cl{R}$ of $\Lambda$, we say that a point $x\in \mathrm{Spec}(\cl{R})(\overline{\bb{Q}}_p)$ is \emph{arithmetic} if it lies above an arithmetic point $P_{r,\epsilon}\in \mathrm{Spec}(\Lambda)(\overline{\bb{Q}}_p)$. In this case, we refer to $k_x=r+2$ as the \emph{weight} of $x$.

Let $N$ be a positive integer coprime to $p$ and let $\chi:(\bb{Z}/Np\bb{Z})^\times \rightarrow \overline{\bb{Q}}^\times$ be a Dirichlet character. A \emph{Hida family} of tame level $N$ and character $\chi$ is a formal $q$-expansion
\[
\bm{\xi}=\sum_{n\geq 1} a_n(\bm{\xi}) q^n \in \Lambda_{\bm{\xi}}[[q]]
\]
with coefficients in a normal domain $\Lambda_{\bm{\xi}}$ finite flat over $\Lambda$ such that for every arithmetic point $x\in \mathrm{Spec}(\Lambda_{\bm{\xi}})(\overline{\bb{Q}}_p)$ lying above some $P_{r,\epsilon}$ with $r\geq 0$ the corresponding specialization $\bm{\xi}_x$ is a $p$-ordinary eigenform in $S_{r+2}(Np^s,\chi \epsilon \omega^{-r})$, where $s=\max\lbrace 1,\ord_p(\mathrm{cond}(\epsilon)\rbrace$. We define $\cl{O}_{\bm{\xi}}=\Lambda_{\bm{\xi}}[1/p]$. We set $U_{\bm{\xi}}=\mathrm{Spec}(\Lambda_{\bm{\xi}})$ and $\cl{U}_{\bm{\xi}}=\mathrm{Spec}(\cl{O}_{\bm{\xi}})$ and we denote by $\kappa_{\bm{\xi}}:U_{\bm{\xi}}\rightarrow \mathrm{Spec}(\Lambda)$ the homomorphism corresponding to the inclusion $\Lambda \xhookrightarrow{} \Lambda_{\bm{\xi}}$. We denote by $\cl{W}_{\bm{\xi}}$ the set of arithmetic points in $U_{\bm{\xi}}(\overline{\bb{Q}}_p)$. We say that an arithmetic point $x\in \cl{W}_{\bm{\xi}}$ of weight $k_x\geq 1$ is \emph{classical} if the corresponding specialization $\bm{\xi}_x$ is a classical modular form (which is always the case if $k_x\geq 2$). In this case $\bm{\xi}_x$ is a $p$-stabilized newform. If $\bm{\xi}_x$ is the $p$-stabilization of a newform of level $N$, we denote the latter by $\bm{\xi}_x^{\circ}$. Otherwise, we write $\bm{\xi}_x^\circ=\bm{\xi}_x$. We note that, if $k_x>2$, then $\bm{\xi}_x$ is always $p$-old.

Attached to a Hida family $\bm{\xi}$ 
there is a locally-free rank-two $\Lambda_{\bm{\xi}}$-module $\bb{V}_{\bm{\xi}}$ equipped with a continuous action of $G_{\bb{Q}}$ such that for any arithmetic point $x\in \cl{W}_{\bm{\xi}}$ of weight $k_x\geq 2$ the specialization $\bb{V}_{\bm{\xi}}\otimes_{\Lambda_{\bm{\xi}},x}\overline{\bb{Q}}_p$ recovers the $G_{\bb{Q}}$-representation $V_{\bm{\xi}_x^\circ}$ attached to $\bm{\xi}_{x}^\circ$. We call $\bb{V}_{\bm{\xi}}$ the \emph{big Galois representation} attached to $\bm{\xi}$. As a $G_{\bb{Q}_p}$-representation, $\bb{V}_{\bm{\xi}}$ admits a filtration
\[
0\rightarrow \bb{V}_{\bm{\xi}}^+\rightarrow \bb{V}_{\bm{\xi}}\rightarrow \bb{V}_{\bm{\xi}}^-\rightarrow 0,
\]
where $\bb{V}_{\bm{\xi}}^{+}$ is a free rank-one $\Lambda_{\bm{\xi}}$-module, $\bb{V}_{\bm{\xi}}^-$ is a locally-free rank-one $\Lambda_{\bm{\xi}}$-module, and the $G_{\bb{Q}_p}$-action on $\bb{V}_{\bm{\xi}}^-$ is unramified with $\Fr_p$ acting as multiplication by $a_p(\bm{\xi})$. For any arithmetic point $x\in \cl{W}_{\bm{\xi}}$ of weight $k_x\geq 2$, the corresponding specialization of the exact sequence above recovers the filtration
\[
0\rightarrow V_{\bm{\xi}_x^\circ}^+\rightarrow V_{\bm{\xi}_x^\circ}\rightarrow V_{\bm{\xi}_x^\circ}^-\rightarrow 0.
\]
Moreover, if $\bm{\xi}$ is residually irreducible and $p$-distinguished, then the $\Lambda_{\bm{\xi}}$-modules $\bb{V}_{\bm{\xi}}$ and $\bb{V}_{\bm{\xi}}^{-}$ are free.

\subsection{The Coleman--Mazur eigencurve}

Let $\bm{\xi}$ be a Hida family of tame level $N_{\bm{\xi}}$ and character $\chi_{\bm{\xi}}$. Let $x\in \cl{W}_{\bm{\xi}}$ be a classical point of weight $k_x\geq 1$. We say that $x$ is a \emph{crystalline} point if $\bm{\xi}_x$ has nebentypus of conductor coprime to $p$. We say that $x$ is a \emph{good} point if in addition any of the following conditions holds: 
\begin{enumerate}
\item[(a)] $k_x \geq 2$;
\item[(b)] $k_x=1$ and $\bm{\xi}_{x}$ is a $p$-stabilisation of a classical $p$-regular cuspidal weight-one newform of level $N_{\bm{\xi}}$ without real multiplication by a real  quadratic field $F$ where the prime $p$ splits;
\item[(c)] $k_x=1$ and $\bm{\xi}_{x}$ is the $p$-stabilisation of a $p$-irregular weight-one Eisenstein series of conductor $N_{\bm{\xi}}$.
\end{enumerate}

Note that these conditions imply in particular that $\bm{\xi}_{x}$ is an \'etale point of the cuspidal Coleman--Mazur--Buzzard $p$-adic eigencurve of tame level $N_{\bm{\xi}}$. When (a) holds, this follows from the work of Hida and Coleman (see, e.g. \cite{bellaiche11a}); when (b) holds, this is proved in \cite{bellaiche-dimitrov}; when (c) holds, it is proved in \cite{BeDiPo}.

If $x\in \cl{W}_{\bm{\xi}}$ is a good point of weight $1$, then it follows from \cite[Prop.~2.2]{BDV} that $\bb{V}_{\bm{\xi}}\otimes_{\Lambda_{\bm{\xi}},x}\overline{\bb{Q}}_p$ recovers the Deligne--Serre Artin representation $V_{\bm{\xi}_x^\circ}$ attached to $\bm{\xi}_{x}^\circ$. In particular, this allows us to obtain a filtration
\[
0\rightarrow V_{\bm{\xi}_x^\circ}^+ \rightarrow V_{\bm{\xi}_x^\circ} \rightarrow V_{\bm{\xi}_x^\circ}^- \rightarrow 0
\]
for $V_{\bm{\xi}_x^\circ}$ as a $G_{\bb{Q}_p}$-representation by specializing the filtration
\[
0\rightarrow \bb{V}_{\bm{\xi}}^+ \rightarrow \bb{V}_{\bm{\xi}} \rightarrow \bb{V}_{\bm{\xi}}^- \rightarrow 0
\]
at the point $x$.

\subsection{Irregular weight-one Eisenstein series}\label{subsec:wt1}

Let $K$ be an imaginary quadratic field of discriminant $-d_K$ in which $p$ splits. Write $p \mathcal O_K = \mathfrak p \bar{\mathfrak p}$, with $\mathfrak p$ the prime of $\mathcal O_K$ corresponding to the fixed embedding $\iota_p \colon \overline{\mathbb Q} \hookrightarrow \overline{\mathbb Q}_p$. Let $\varepsilon_K$ stand for the quadratic character attached to the imaginary quadratic field $K$. In particular, for any prime $q\nmid d_K$, we have $\varepsilon_K(q) = +1$ if $q$ splits in $K$ and $\varepsilon_K(q) = -1$ if $q$ is inert in $K$.

\begin{defi}
The weight-one Eisenstein series of character $\varepsilon_K$ is the weight-one modular form given by the Fourier expansion \[ \Eis_1(\varepsilon_K) = \frac{1}{2} L(\varepsilon_K,0) + \sum_{n \geq 1} q^n \sum_{d|n} \varepsilon_K(d) \in M_1(-d_K, \varepsilon_K). \] It has level $\Gamma_1(-d_K)$ and character $\varepsilon_K$. 
\end{defi} 

We will write $f=\Eis_1(\varepsilon_K)$. Note that the eigenform $f$ is a $p$-irregular modular form, i.e., its $p$-th Hecke polynomial, given by $X^2-2X+1$, has a double root. Define \[ f_\alpha = f(q) - f(q^p) \in M_1(-pd_K, \varepsilon_K) \] as the unique $p$-stabilisation of $f$.

A result of Betina, Dimitrov and Pozzi \cite{pozzi-thesis}, \cite{BeDiPo} extending the prior work of Bella\"iche and Dimitrov \cite{bellaiche-dimitrov}, establishes the following.

\begin{propo}
The eigenform $f_\alpha$ is an \'etale point of the cuspidal Coleman--Mazur eigencurve. In particular, there exists a unique (up to conjugation) Hida family $\hf$ of tame level $-d_K$ and tame character $\chi_f = \varepsilon_K$ passing through $f_\alpha$. Moreover, the Hida family $\hf$ has complex multiplication by $K$.
\end{propo}
\begin{proof}
This is \cite[Thm. A]{BeDiPo}.
\end{proof}

\begin{remark}
Similar results exist if we consider instead the Eisenstein series $\Eis_1(\chi_1,\chi_2 \varepsilon_K)$, with $\chi_1$ and $\chi_2$ Dirichlet characters of prime-to-$p$ conductors $N_1$ and $N_2$ with $\chi_1(p)=\chi_2(p)$. In that case, one would need to define $f_\alpha$ as \[ f_\alpha = \Eis_1(\varepsilon_K)(q) - \chi_1(p) \Eis_1(\varepsilon_K)(q^p) \in M_1(-pd_KN_1N_2, \varepsilon_K \chi_1 \chi_2). \]
\end{remark}


Let $K^{\frk{p}}_{\infty}$ be the $\bb{Z}_p$-extension of $K$ unramified away from $\frk{p}$. Let $\Gamma_{\frk{p}}=\Gal(K^{\frk{p}}_{\infty}/K)$ and let $\Lambda_\hf=\cl{O}[[\Gamma_{\frk{p}}]]$. Let $\rec_K:\bb{A}_K^\times\rightarrow G_K^{\ab}$ denote the (arithmetically normalized) global reciprocity map and let $\Theta_\frk{p}$ denote the composition of $\rec_K$ with the canonical projection $G_K^{\ab} \twoheadrightarrow \Gamma_{\frk{p}}$. Let $\rec_{K,\frk{p}}:K_{\frk{p}}^\times \rightarrow G_{K_{\frk{p}}}^{\ab}$ denote the (arithmetically normalized) local reciprocity map. We identify $G_{K_\frk{p}}$ with the decomposition group above $\frk{p}$ determined by our fixed embedding $\overline{\bb{Q}}\hookrightarrow \overline{\bb{Q}}_p$, and thus we regard $G_{K_\frk{p}}$ as a subgroup of $G_K$. Thus we obtain a homomorphism $G_{K_\frk{p}}^{\ab} \rightarrow G_K^{\ab} \twoheadrightarrow \Gamma_{\frk{p}}$ (which is surjective if $p\nmid h_K$) and we define $\theta_\frk{p}:K_\frk{p}^\times\rightarrow \Gamma_{\frk{p}}$ as the composition of $\rec_{K,\frk{p}}$ with this homomorphism. The character $1+p\bb{Z}_p\rightarrow \Lambda_{\hf}^\times$ defined by $z\mapsto z^{-1}\theta_{\frk{p}}(z)$ extends to an embedding $\Lambda\hookrightarrow \Lambda_{\hf}$, whereby $\Lambda_{\hf}$ becomes a finite flat extension of $\Lambda$. For each non-zero fractional ideal $\frk{a}$ of $K$, let $x_\frk{a}\in \bb{A}_{K,f}^\times$ be a finite id\`ele with $\ord_w(x_{\frk{a},w})=\ord_w(\frk{a})$ at each finite place $w$ with $\ord_w(\frk{a})\neq 0$ and $x_{\frk{a},w}=1$ for all other places $w$. Then, the Hida family $\hf$ is given by
\[
\hf=\sum_{\frk{a}\subseteq \cl{O}_K,\,\frk{p}\nmid\frk{a}}\Theta_\frk{p}(x_\frk{a}) q^{N_{K/\bb{Q}}(\frk{a})} \in \Lambda_\hf[[q]].
\]
The specialization of $\hf$ at the point $x_0\in \cl{W}_{\hf}$ corresponding to the trivial character of $\Gamma_{\frk{p}}$ recovers the modular form $f_\alpha$.

Let $\bm{\varphi}:G_K\twoheadrightarrow \Gamma_{\frk{p}}$ denote the canonical projection. Then, we have the following result.

\begin{proposition}\label{prop:BSTW}
The $\Lambda_{\hf}[G_{\mathbb Q}]$-modules $\bb{V}_{\hf}$ and $\Ind_K^{\mathbb Q} \Lambda_\hf(\bm{\varphi})$ are isomorphic.
\end{proposition}
\begin{proof}
This is \cite[Thm.~1.21]{BSTW}. 
\end{proof}

Fix once and for all an isomorphism of $\Lambda_{\hf}[G_{\mathbb Q}]$-modules 
\begin{equation}\label{eq:iso}
\gamma \, : \, \bb{V}_{\hf} \cong \Ind_K^{\mathbb Q} \bm{\varphi}.
\end{equation}
This isomorphism is not canonical, so the next decompositions will depend on this choice.
Since $p$ splits in $K$, the restrictions of $\Ind_K^{\mathbb Q} \bm{\varphi}$ to $G_K$ and $G_{\mathbb Q_p}$ decompose as the direct sum of $\bm{\varphi}$ and its complex conjugate $\bm{\varphi}^c$. Note that the character $\bm{\varphi}^c|_{G_{\mathbb Q_p}}$ is unramified and maps the arithmetic Frobenius $\Fr_p$ to the $p$-th Fourier coefficient $a_p(\hf) = \Theta(x_{\bar{\mathfrak p}})$ of $\hf$. Therefore, the restriction of $\bb{V}_{\hf}$ to $G_{\bb{Q}_p}$ decomposes
\[
\bb{V}_{\hf} = \bb{V}_{\hf}^+ \oplus \bb{V}_{\hf}^-, \quad \text{ with } \quad \gamma(\bb{V}_{\hf}^+) = \Lambda_{\hf}(\bm{\varphi}|_{G_{\mathbb Q_p}}) \text{ and } \gamma(\bb{V}_{\hf}^-) = \Lambda_\hf(\bm{\varphi}^c|_{G_{\mathbb Q_p}}).
\]
Recall that the specialization $\bb{V}_{\hf}\otimes_{\Lambda_{\hf},x_0} L$ recovers the Deligne--Serre Artin representation $V_f$. Setting, as before, $V_f^{\pm}=\bb{V}_{\hf}^{\pm}\otimes_{\Lambda_{\hf},x_0} L$, we have a decomposition $V_f=V_f^+\oplus V_f^-$ of $G_{\bb{Q}_p}$-representations. Also, specializing the isomorphism $\gamma$ at $x_0$, we obtain an isomorphism
\[
\gamma:V_f\cong L(\bm{1})\oplus L(\varepsilon_K) 
\]
of $L[G_{\bb{Q}}]$-modules.

Let $\bm{v}^+$ and $\bm{v}^-$ be the canonical $\Lambda_{\hf}$-bases of the $G_K$-submodules $\Lambda_{\hf}(\bm{\varphi})$ and $\Lambda_{\hf}(\bm{\varphi}^c)$ of $\Ind_K^{\mathbb Q} \Lambda_\hf(\bm{\varphi})$, respectively. The maps $\bm{v}^{\pm} \colon G_{\mathbb Q} \rightarrow \Lambda_{\hf}$ are determined by $(\bm{v}^+(1),\bm{v}^+(c))=(1,0)$ and $(\bm{v}^-(1),\bm{v}^-(c)) = (0,1)$, where $c$ denotes our fixed choice of complex conjugation. Set $\bm{v}_\hf^{\pm} = \gamma^{-1}(\bm{v}^{\pm})$ in $\bb{V}_{\hf}^{\pm}$, let $v_f^{\pm}$ in $V_f^{\pm}$ be their specializations at $x_0$ and define \[ v_{f,1} = v_f^+ + v_f^-, \quad v_{f,\varepsilon_K} = v_f^+ - v_f^-. \] Note that $c$ exchanges the vectors $\bm{v}^+$ and $\bm{v}^-$, and therefore the elements $\gamma(v_{f,1})$ and $\gamma(v_{f,\varepsilon_K})$ give bases of the $G_{\mathbb Q}$-representations $L(\bm{1})$ and $L(\varepsilon_K)$, respectively. These choices of vectors will be used later on.

\section{$p$-adic $L$-functions}

Let $(\hg,\hh)$ be a pair of cuspidal Hida families and let $\hf$ be the Hida family introduced in \S\ref{subsec:wt1}. This work relies on the comparison between the Hida--Rankin $p$-adic $L$-function $L_p(\hg,\hh)$ and the triple product $p$-adic $L$-function $L_p(\hf,\hg,\hh)$. In this section we introduce these objects.



\subsection{Hida--Rankin $p$-adic $L$-function}\label{subsec:Hida-Rankin}


Here we introduce the three-variable Hida--Rankin $p$-adic $L$-function. We follow mainly the exposition in \cite{KLZ}.

Let $(\hg,\hh)$ be a pair of cuspidal Hida families of tame levels $(N_g,N_h)$ and characters $(\chi_{\hg},\chi_{\hh})$. We make the following assumption.

\begin{ass}
    $\gcd(N_g,N_h)=1$.
\end{ass}

Put $\cl{W}_{\hg \hh \mathbf{s}}\coloneqq \cl{W}_{\hg}\times \cl{W}_{\hh}\times \cl{W}$.
For $\phi\in\lbrace \hg,\hh\rbrace$, let $\cl{W}_{\phi}^\circ$ denote the subset of crystalline points in $\cl{W}_{\phi}$. Also, fix an integer $s_0$ and define
\[
\cl{W}^\circ=\lbrace P_{s,\bm{1}}\in \cl{W}\;\vert\; s\equiv s_0 \pmod{p-1}\rbrace.
\]
In an abuse of notation, a point $P_{s,\bm{1}}\in \cl{W}^\circ$ might be denoted simply by $s$.

Let $(y_0,z_0)\in \cl{W}_{\hg}^\circ\times \cl{W}_{\hh}^\circ$ be a good crystalline point with corresponding weights $(l_0,m_0)$ satisfying $l_0\geq 2$, $m_0\geq 1$. We define
\begin{align*}
\cl{W}_{\hg}^{\mathrm{cl}}&\coloneqq \lbrace y\in \cl{W}_{\hg}^\circ \;\vert\; k_y\geq 2 \rbrace; \\
\cl{W}_{\hh}^{\mathrm{cl}}&\coloneqq \lbrace z\in \cl{W}_{\hh}^\circ \;\vert\; k_z\geq 2\rbrace \cup \lbrace z_0\rbrace.
\end{align*}
Note that, for $\phi\in\lbrace \hg,\hh\rbrace$, the specialization of $\phi$ at a point in $\cl{W}_{\phi}^{\mathrm{cl}}$ is a classical $p$-stabilized newform. For each $x\in \cl{W}_{\phi}^{\mathrm{cl}}$, we define $\alpha_{\phi_x}\coloneqq a_p(\phi_x)$ and $\beta_{\phi_x}\coloneqq\chi_{\phi_x}(p)p^{k_x-1}\alpha_{\phi_x}^{-1}$. Note that, if $k_x\geq 2$, then $\alpha_{\phi_x}$ is the unit root of the $p$-th Hecke polynomial of $\phi_x^{\circ}$.

Let $\cl{W}_{\hg \hh \mathbf{s}}^{\mathrm{cl}}=\cl{W}_{\hg}^{\mathrm{cl}}\times \cl{W}_{\hh}^{\mathrm{cl}}\times \cl{W}^\circ$ and let $\cl{W}_{\hg \hh \mathbf{s}}^{g}$ denote the subset of points $(y,z,s)\in \cl{W}_{\hg \hh \mathbf{s}}^{\mathrm{cl}}$ of weights $(k_y,k_z)=(l,m)$ satisfying $m\leq s<l$. This is the range of interpolation for the three-variable Rankin $p$-adic $L$-function $L_p(\hg,\hh,\mathbf{s})$ introduced below. Similarly, we define $\cl{W}_{\hg\hh\mathbf{s}}^h$ as the subset of points $(y,z,s)\in \cl{W}_{\hg \hh \mathbf{s}}^{\mathrm{cl}}$ of weights $(k_y,k_z)=(l,m)$ satisfying $l\leq s<m$. Note that both $\cl{W}_{\hg\hh\mathbf{s}}^g$ and $\cl{W}_{\hg\hh\mathbf{s}}^h$ exclude points with $l=m$. Set $\Lambda_{\hg \hh \mathbf{s}}=\Lambda_{\hg}\hat{\otimes}_{\cl{O}} \Lambda_{\hh}\hat{\otimes}_{\bb{Z}_p} \Lambda$ and $\cl{O}_{\hg \hh \mathbf{s}}=\Lambda_{\mathbf{hgs}}[1/p]$. Let $I_\hg\subseteq \Lambda_\hg$ be the congruence ideal of $\hg$ and let $I_{\hh}\subseteq \Lambda_{\hh}$ be the congruence ideal of $\hh$.




 Hida constructed in \cite{hida1} and \cite{hida2} a three-variable $p$-adic Rankin $L$-function $L_p(\hg, \hh, \mathbf{s})$ in $I_{\hg}^{-1} \cl{O}_{\mathbf{ghs}}$ interpolating the algebraic parts of the critical values $L(\hg^\circ_y,\hh^\circ_z,s)$ for every triple of critical points $(y,z,s)$ in $\mathcal W_{\hg \hh \mathbf{s}}^{g}$. The following formulation is taken from \cite[Thm.~7.7.2]{KLZ}. 

\begin{theorem}(Hida)\label{thm:Hida-3var}
There exists a $p$-adic $L$-function
\[
L_p(\hg,\hh,\mathbf{s}) \in I_{\hg}^{-1}\cl{O}_{\hg \hh \mathbf{s}}
\]
such that for every $(y,z,s)\in\mathcal W_{\hg \hh}^g$ of weights $(k_y,k_z)=(l,m)$ with $p\nmid \mathrm{cond}(\hg_y)\cdot \mathrm{cond}(\hh_z)$ we have
\[
L_p(\hg,\hh,\mathbf{s})(y,z,s) = \frac{\mathcal{E}(y,z,s)}{\mathcal{E}_0(y)\mathcal{E}_1(y)} \frac{\Gamma(s)\Gamma(s-m+1)}{\pi^{2s-m+1}(-i)^{l-m}2^{2s+l-m} \langle \hg_y^\circ,\hg_y^\circ \rangle}  \times L(\hg_y^\circ,\hh_z^\circ,s)
\]
where $\langle \hg_y^{\circ}, \hg_y^{\circ} \rangle$ is the Petersson norm as normalized in \cite{KLZ} and the Euler factors are defined by
\begin{eqnarray}
      \nonumber \mathcal E_0(y) &:=& 1-\chi_g^{-1}(p) \beta_{\hg_y}^2p^{1-l}, \\
      \nonumber \mathcal E_1(y) &:=& 1-\chi_g(p)\alpha_{\hg_y}^{-2}p^{l-2}, \\
      \mathcal E(y,z,s) &:=& \Big(1-\frac{p^{s-1}}{\alpha_{\hg_y} \alpha_{\hh_z}} \Big) \Big(1-\frac{p^{s-1}}{\alpha_{\hg_y} \beta_{\hh_z}} \Big) \Big(1-\frac{\beta_{\hg_y} \alpha_{\hh_z}}{p^{s}}\Big) \Big(1-\frac{\beta_{\hg_y} \beta_{\hh_z}}{p^{s}}\Big).
\end{eqnarray}
\end{theorem}

\begin{remark}
Similarly, there is a $p$-adic $L$-function $L_p(\hh,\hg,\mathbf{s}) \in I_{\hh}^{-1}\cl{O}_{\hg \hh \mathbf{s}}$ interpolating the Rankin--Selberg $L$-values in the region $\cl{W}_{\hg\hh\mathbf{s}}^h$. Note that the $p$-adic $L$-functions $L_p(\hg,\hh,\mathbf{s})$ and $L_p(\hh,\hg,\mathbf{s})$ are different. 
\end{remark}

This construction was later generalized by Urban to the case where the modular forms are not necessarily ordinary \cite{Urban-nearly-overconvergent}, \cite[App.]{AI}, using nearly overconvergent modular forms of finite order and their spectral theory. In this case, the interpolation property does not completely determine the $p$-adic $L$-function, and one needs to impose further conditions.

\subsection{Triple product $p$-adic $L$-function}\label{subsec:tripleproduct}

Let $(\hf,\hg,\hh)$ be a triple of Hida families of tame levels $(N_f,N_g,N_h)$ and characters $(\chi_{\hf},\chi_{\hg},\chi_{\hh})$. In this section, we introduce the three-variable $p$-adic $L$-function interpolating the central values of the triple product $L$-functions associated with classical specializations of $(\hf,\hg,\hh)$ in the $g$-dominant region. 

Put $\cl{W}_{\hf \hg \hh}\coloneqq \cl{W}_{\hf}\times \cl{W}_{\hg}\times \cl{W}_{\hh}$. For $\phi\in\lbrace \hf,\hg,\hh\rbrace$, let $\cl{W}_{\phi}^\circ$ denote the subset of crystalline points in $\cl{W}_{\phi}$ and put $\cl{W}_{\hf \hg \hh}^{\circ}=\cl{W}_{\hf}^{\circ} \times \cl{W}_{\hg}^{\circ} \times \cl{W}_{\hh}^{\circ}\subseteq \cl{W}_{\hf \hg \hh}$. Let $w_0=(x_0,y_0,z_0)\in \cl{W}_{\hf \hg \hh}^\circ$ be a good crystalline point with corresponding weights $(k_0,l_0,m_0)$ satisfying $k_{0}+l_{0}+m_{0}\equiv 0\pmod{2}$, $k_0,m_0\geq 1$, $l_0\geq 2$. Put $(f,g,h)\coloneqq (\hf_{x_0}^\circ,\hg_{y_0}^\circ,\hh_{z_0}^\circ)$. Note that $(f,g,h)$ is a triple of $p$-stabilized newforms with characters $(\chi_f,\chi_g,\chi_h)=(\chi_{\hf}^{(p)},\chi_{\hg}^{(p)},\chi_{\hh}^{(p)})$, where $\chi^{(p)}$ stands for the prime-to-$p$ part of the character $\chi$. 

We make the following assumption.

\begin{ass} 
\hfill
    \begin{enumerate}
        \item $\gcd(N_f,N_g,N_h)$ is square-free;
        \item $\chi_f\chi_g\chi_h=1$;
        \item the residual representation $\bar{\rho}_{g}$ is absolutely irreducible and $p$-distinguished.
    \end{enumerate}
\end{ass}

Put $V_{fgh}\coloneqq V_f\otimes_L V_g\otimes_L V_h$. Then $V_{fgh}$ is a $G_\bb{Q}$-representation of dimension $8$ over $L$. Moreover, putting $c_0\coloneqq(k_{0}+l_{0}+m_0-2)/2$, the twisted representation $V_{fgh}^\dagger\coloneqq V_{fgh}(1-c_0)$ is Kummer self-dual, i.e., $(V_{fgh}^\dagger)^\vee\simeq V_{fgh}^\dagger(1)$.


We now add the following assumption.

\begin{ass}
    For all prime $\ell\mid N_fN_gN_h$, $\epsilon_\ell(V_{fgh}^\dagger)=+1$.
\end{ass}

We define
\begin{align*}
\tilde{\cl{W}}_{\hf}^{\mathrm{cl}}&\coloneqq \lbrace x\in \cl{W}_{\hf}^\circ \;\vert\; k_x\geq 2 \text{ and } k_x\equiv k_0\;(\mathrm{mod}\,2(p-1))\rbrace \cup \lbrace x_0\rbrace; \\ 
\tilde{\cl{W}}_{\hg}^{\mathrm{cl}}&\coloneqq \lbrace y\in \cl{W}_{\hg}^\circ \;\vert\; k_y\geq 2 \text{ and } k_y\equiv l_0 \;(\mathrm{mod}\,2(p-1))\rbrace; \\
\tilde{\cl{W}}_{\hh}^{\mathrm{cl}}&\coloneqq \lbrace z\in \cl{W}_{\hh}^\circ \;\vert\; k_z\geq 2 \text{ and } k_z\equiv m_0 \;(\mathrm{mod}\,2(p-1))\rbrace \cup \lbrace z_0\rbrace.
\end{align*}
Note that, for $\phi\in\lbrace \hf,\hg,\hh\rbrace$, the specialization of $\phi$ at a point in $\tilde{\cl{W}}_{\phi}^{\mathrm{cl}}$ is a classical $p$-stabilized newform. For each $x\in \tilde{\cl{W}}_{\phi}^{\mathrm{cl}}$, we define $\alpha_{\phi_x}\coloneqq a_p(\phi_x)$ and $\beta_{\phi_x}\coloneqq\chi_{\phi_x}(p)p^{k_x-1}\alpha_{\phi_x}^{-1}$. Note that, if $k_x\geq 2$, then $\alpha_{\phi_x}$ is the unit root of the $p$-th Hecke polynomial of $\phi_x^{\circ}$.

Put $\cl{W}_{\hf \hg \hh}^{\mathrm{cl}}=\tilde{\cl{W}_{\hf}}^{\mathrm{cl}}\times \tilde{\cl{W}}_{\hg}^{\mathrm{cl}}\times \tilde{\cl{W}}_{\hh}^{\mathrm{cl}}$. This set admits the natural partition
\[
\mathcal W_{{\bf fgh}}^{\mathrm{cl}} = \mathcal W_{{\bf fgh}}^{f} \sqcup \mathcal W_{{\bf fgh}}^{g} \sqcup \mathcal W_{{\bf fgh}}^{h} \sqcup \mathcal W_{{\bf fgh}}^{\bal},
\]
where
\begin{itemize}
\item $\mathcal W_{{\bf fgh}}^g$ denotes the set of points $(x,y,z) \in \mathcal W_{{\bf fgh}}^{\mathrm{cl}}$ of weights $(k,l,m)$ such that $l \geq k+m$;
\item $\mathcal W_{{\bf fgh}}^f$ (resp. $\mathcal W_{{\bf fgh}}^h$) is defined similarly, replacing the role of $\hg$ by $\hf$ (resp. $\hh$);
\item $\mathcal W_{{\bf fgh}}^{\bal}$ is the set of balanced triples, consisting of points $(x,y,z)\in \mathcal W_{{\bf fgh}}^{\mathrm{cl}}$ of weights $(k,l,m)$ such that each of the weights is strictly smaller than the sum of the other two.
\end{itemize}
Set $\Lambda_{\hf \hg \hh}=\Lambda_{\hf}\hat{\otimes}_{\cl{O}} \Lambda_{\hg}\hat{\otimes}_{\cl{O}} \Lambda_{\hh}$ and $\cl{O}_{\hf \hg \hh}=\Lambda_{\hf \hg \hh}[1/p]$. Let $I_\hg\subseteq \Lambda_\hg$ be the congruence ideal of $\hg$.


In \cite{DR1}, Darmon and Rotger constructed a triple product $p$-adic $L$-function interpolating the square roots of the central values of the triple product $L$-functions associated with classical specializations of $(\hf,\hg,\hh)$ in $\mathcal W_{{\bf fgh}}^g$. The precise interpolation formula given below is due to Hsieh \cite{Hs}

\begin{theorem}(Darmon--Rotger, Hsieh)\label{thm:hsieh}
There exists a $p$-adic $L$-function
\[
\Lp^g(\hf,\hg,\hh) \in I_\hg^{-1}\cl{O}_{\hf \hg \hh}
\]
such that for every $(x,y,z) \in \mathcal W_{{\bf fgh}}^g$ of weights $(k,l,m)$ with $p\nmid \mathrm{cond}(\hf_x)\cdot\mathrm{cond}(\hg_y)\cdot\mathrm{cond}(\hh_z)$ we have
\[
\Lp^g(\hf,\hg,\hh)^2(x,y,z) = \frac{\mathfrak a(k,l,m)}{\langle \hg_y^{\circ},\hg_y^{\circ} \rangle^2} \cdot \mathfrak e^2(x,y,z) \cdot L(\hf_x^{\circ},\hg_y^{\circ},\hh_z^{\circ},c),
\]
where
\begin{enumerate}
\item $c = \frac{k+l+m-2}{2}$,
\item $\mathfrak a(k,l,m) = (2\pi i)^{-2(l-2)} \cdot \Big( \frac{k+l+m-4}{2} \Big) ! \cdot \Big( \frac{-k+l+m-2}{2} \Big) ! \cdot \Big( \frac{k+l-m-2}{2} \Big) ! \cdot \Big( \frac{-k+l-m}{2} \Big) !$,
\item $\mathfrak e(x,y,z) = \mathcal E(x,y,z)/\mathcal E_0(x) \mathcal E_1(x)$ with
    \begin{eqnarray}
      \nonumber \mathcal E_0(y) &:=& 1-\chi_g^{-1}(p) \beta_{\hg_y}^2p^{1-l}, \\
      \nonumber \mathcal E_1(y) &:=& 1-\chi_g(p)\alpha_{\hg_y}^{-2}p^{l-2}, \\
      \nonumber \mathcal E(x,y,z) &:=& \Big(1-\chi_g(p) \alpha_{\hf_x} \alpha_{\hg_y}^{-1}\alpha_{\hh_z}p^{\frac{-k+l-m}{2}} \Big) \times \Big(1-\chi_g(p) \alpha_{\hf_x} \alpha_{\hg_y}^{-1}\beta_{\hh_z}p^{\frac{-k+l-m}{2}} \Big)\\ \nonumber
       & & \times \Big(1-\chi_g(p) \beta_{\hf_x} \alpha_{\hg_y}^{-1}\alpha_{\hh_z}p^{\frac{-k+l-m}{2}} \Big) \times \Big(1-\chi_g(p) \beta_{\hf_x} \alpha_{\hg_y}^{-1}\beta_{\hh_z}p^{\frac{-k+l-m}{2}} \Big).
    \end{eqnarray}
\end{enumerate}
\end{theorem}

\begin{remark}
    To be precise, each choice of test vectors $(\breve{\hf},\breve{\hg},\breve{\hh})$ for $(\hf,\hg,\hh)$ determines a $p$-adic $L$-function $\Lp^g(\breve{\hf},\breve{\hg},\breve{\hh})$, and Hsieh shows that there exists an optimal choice of test vectors $(\breve{\hf}^\ast,\breve{\hg}^\ast,\breve{\hh}^\ast)$ for which the $p$-adic $L$-function $\Lp^g(\breve{\hf}^\ast,\breve{\hg}^\ast,\breve{\hh}^\ast)$, which we denote simply by $\Lp^g(\hf,\hg,\hh)$, satisfies the precise interpolation formula given above. Using these same test vectors, one can also define $p$-adic $L$-functions $\Lp^f(\hf,\hg,\hh)$ and $\Lp^h(\hf,\hg,\hh)$ interpolating square roots of classical $L$-values in the regions $\mathcal W_{{\bf fgh}}^f$ and $\mathcal W_{{\bf fgh}}^h$, respectively. Note that our choice of test vectors might not be optimal for these regions, so $\Lp^f(\hf,\hg,\hh)$ and $\Lp^h(\hf,\hg,\hh)$ might not precisely satisfy an interpolation formula analogous to the one above. However, if $\gcd(N_g,N_h)=1$, $\hg$ and $\hh$ are residually irreducible and $p$-distinguished, and $\hf$ is the CM Hida family introduced in \S\ref{subsec:wt1}, which is the situation that we will consider, then it follows from \cite{Hs} that one can choose test vectors which are optimal for both $\mathcal W_{{\bf fgh}}^g$ and $\mathcal W_{{\bf fgh}}^h$.
\end{remark}




\section{Families of cohomology classes}

In this section we introduce the two kinds of cohomology classes that are involved in our comparison, together with the corresponding reciprocity laws connecting them with the $p$-adic $L$-functions introduced in the previous section.

The cohomology classes introduced in this section are the following:
\begin{enumerate}
\item[(a)] the Beilinson--Flach classes attached to a pair of Hida families $(\hg,\hh)$ constructed in \cite{KLZ};
\item[(b)] the diagonal cycles attached to a triple of Hida families $(\hf,\hg,\hh)$ constructed in \cite{BSV} and \cite{DR3}. 
\end{enumerate}


\subsection{Beilinson--Flach classes}\label{subsec:BF}

We use the notations and assumptions introduced in \S\ref{subsec:Hida-Rankin}. In particular, $(\hg,\hh)$ is a pair of Hida families of tame levels $(N_g,N_h)$ and characters $(\chi_\hg,\chi_\hh)$. 
We also make the following assumption.

\begin{ass}
    The Hida families $\hg$ and $\hh$ are residually irreducible and $p$-distinguished.
\end{ass}

Let $\mathbf{s}:\bb{Z}_p^\times\rightarrow \Lambda^\times$ be the character defined by $z\mapsto \omega^{s_0}(z)[\langle z\rangle]$, where $\langle z\rangle =\omega^{-1}(z)z$ and $[\langle z\rangle]$ denotes the corresponding group-like element. Note that the specialization of $\mathbf{s}$ at a point $P_{s,\bm{1}}\in \cl{W}^\circ$ is the character $z\mapsto z^s$. We denote by $\mathbf{l}:\bb{Z}_p^\times\rightarrow \Lambda_\hg^\times$ the weight character of $\hg$, defined by $z\mapsto \chi_{\hg,p}(z) z^2[\langle z\rangle]_{\Lambda_{\hg}}$. Note that the specialization of $\mathbf{l}$ at a crystalline point $y\in \cl{W}_{\hg}^\circ$ of weight $l$ is the character $z\mapsto z^l$. Similarly, we denote by $\mathbf{m}:\bb{Z}_p^\times \rightarrow \Lambda_\hh^\times$ the weight character of $\hh$, defined by $z\mapsto \chi_{\hh,p}(z)z^2[\langle z\rangle]_{\Lambda_{\hh}}$.

Given specializations $\hg_y$ and $\hh_z$, we denote by $\eta_{\hg_y}^\alpha$ and $\omega_{\hh_z}$ the associated differentials as defined in \cite[\S3]{KLZ}. As shown in \cite[\S10]{KLZ}, 
one can define objects $\bm{\eta}_{\hg}$ and $\bm{\omega}_{\hh}$ interpolating these differentials.


We first introduce a Perrin-Riou logarithm that will be used for the formulation of an explicit reciprocity law. 


\begin{propo}\label{perrin1}
There exists a
homomorphism of $\Lambda_{\hg \hh \mathbf{s}}$-modules
\[
\mathcal L_{\hg \hh \mathbf{s}}^{-+}: H^1(\mathbb Q_p, \mathbb V_{\hg}^- \hat \otimes \mathbb V_{\hh}^+ \hat \otimes \Lambda(1-\mathbf{s})) \rightarrow  I_{\hg}^{-1} I^{-1}\cl{O}_{\hg \hh \mathbf{s}},
\]
where $I$ is the augmentation ideal of $\Lambda$, such that for every point $(y,z,s) \in \mathcal W_{{\bf ghs}}^{\mathrm{cl}}$ of weights $(k_y,k_z)=(l,m)$ with $\alpha_{\hg_y}\beta_{\hh_z}\neq p^s$ the specialization of $\mathcal L_{\hg \hh \mathbf{s}}^{-+}$ at $(y,z,s)$ is the homomorphism
\[
\mathcal L_{{\bf ghs}}^{-+}(y,z,s): H^1(\mathbb Q_p, V_{\hg_y}^- \otimes V_{\hh_z}^+(1-s)) \rightarrow \mathbb C_p
\]
given by
\[
\mathcal L_{{\bf ghs}}^{-+}(y,z,s) = \frac{1-p^{s-1} \alpha_{\hg_y}^{-1} \beta_{\hh_z}^{-1}}{1-p^{-s} \alpha_{\hg_y} \beta_{\hh_z}} \times \begin{cases} \frac{(-1)^{m-s-1}}{(m-s-1)!} \times \langle \log_{\BK}(\cdot), \eta_{\hg_y}^\alpha \otimes \omega_{\hh_z} \rangle & \text{ if } s < m \\ (s-m)! \times \langle \exp_{\BK}^*(\cdot), \eta_{\hg_y}^\alpha \otimes \omega_{\hh_z} \rangle & \text{ if } s \geq m, \end{cases}
\]
where $\log_{\BK}$ is the Bloch--Kato logarithm and $\exp_{\BK}^*$ is the dual exponential map.
\end{propo}
\begin{proof}
This follows from \cite[Theorem~8.2.8, Proposition~10.1.1]{KLZ}.
\end{proof}






\begin{remark}
    We recall that the map in the previous theorem is obtained as the composition
    \[
    \mathcal L_{\hg \hh \mathbf{s}}^{-+} = \langle \log_{\hg \hh \mathbf{s}}^{-+}(\cdot), \bm{\eta}_{\hg} \otimes \bm{\omega}_{\hh} \rangle,
    \]
    where $\log_{\hg \hh \mathbf{s}}^{-+}(\cdot)$ is the Perrin--Riou big logarithm introduced in \cite[Theorem~8.2.8]{KLZ}. 
\end{remark}

\begin{remark}
    Interchanging the roles of $\hg$ and $\hh$, we also have a map
    \[
    \mathcal L_{\hg \hh\mathbf{s}}^{+-}: H^1(\mathbb Q_p, \mathbb V_{\hg}^+ \hat \otimes \mathbb V_{\hh}^- \hat \otimes \Lambda(1-\mathbf{s})) \rightarrow  I_{\hh}^{-1} I^{-1}\cl{O}_{\hg \hh \mathbf{s}},
    \]
    with analogous interpolation properties.
\end{remark}

We also recall the following result.

\begin{propo}\label{prop:local-property}
    The inclusion $\bb{V}_{\hh}^+\xhookrightarrow{} \bb{V}_{\hh}$ induces an injection
    \[
    H^1(\bb{Q}_p,\bb{V}_{\hg}^-\hat{\otimes}\bb{V}_{\hh}^+\hat{\otimes}\Lambda(1-\mathbf{s}))\xhookrightarrow{}H^1(\bb{Q}_p,\bb{V}_{\hg}^-\hat{\otimes}\bb{V}_{\hh}\hat{\otimes}\Lambda(1-\mathbf{s}))
    \]
\end{propo}
\begin{proof}
    This is part of \cite[Prop.~8.1.7]{KLZ}.
\end{proof}

\begin{remark}
    By the previous proposition, we can regard the module $H^1(\bb{Q}_p,\bb{V}_{\hg}^-\hat{\otimes}\bb{V}_{\hh}^+\hat{\otimes}\Lambda(1-\mathbf{s}))$ as a submodule of $H^1(\bb{Q}_p,\bb{V}_{\hg}^-\hat{\otimes}\bb{V}_{\hh}\hat{\otimes}\Lambda(1-\mathbf{s}))$. Similarly, we can also regard the module $H^1(\bb{Q}_p,\bb{V}_{\hg}^+\hat{\otimes}\bb{V}_{\hh}^-\hat{\otimes}\Lambda(1-\mathbf{s}))$ as a submodule of $H^1(\bb{Q}_p,\bb{V}_{\hg}\hat{\otimes}\bb{V}_{\hh}^-\hat{\otimes}\Lambda(1-\mathbf{s}))$.
\end{remark}


To shorten notation, put $\bb{V}_{\hg\hh\mathbf{s}}=\bb{V}_{\hg}\hat{\otimes}\bb{V}_{\hh}\hat{\otimes}\Lambda(1-\mathbf{s})$. Let $\mathscr{F}^2\mathbb V_{\hg\hh\mathbf{s}}\subset \mathbb V_{\hg\hh\mathbf{s}}$ be the $G_{\bb{Q}_p}$-stable $\Lambda_{\hg \hh \mathbf{s}}$-submodule of rank $3$ defined by
\[
\mathscr{F}^2\mathbb V_{\hg\hh\mathbf{s}}=(\bb{V}_{\hg}^+\hat{\otimes}\bb{V}_{\hh} +\bb{V}_{\hg}\hat{\otimes}\bb{V}_{\hh}^+)\hat{\otimes}\Lambda(1-\mathbf{s}).
\]
Fix a finite set $\Sigma$ of places of $\mathbb Q$ containing $\infty$ and the primes dividing $N_gN_hp$ and let $\mathbb Q^\Sigma$ be the maximal extension of $\mathbb Q$ unramified outside $\Sigma$.

\begin{definition}\label{def:Sel1}
The \textit{balanced Selmer group of $\bb{V}_{\hg \hh \mathbf{s}}$} is defined by
\[
H^1_{\bal}(\bb{Q},\mathbb V_{\hg \hh \mathbf{s}}) = \ker \bigg(
H^1(\mathbb Q^\Sigma/\mathbb Q, \mathbb V_{\hg \hh \mathbf{s}}) \longrightarrow
\frac{H^1(\mathbb Q_p,\mathbb V_{\hg \hh \mathbf{s}})}{H^1_{\bal}(\bb{Q}_p,\mathbb V_{\hg \hh \mathbf{s}})}\bigg),
\]
where
\[
H^1_{\bal}(\bb{Q}_p,\mathbb V_{\hg \hh \mathbf{s}})=\ker \big(H^1(\mathbb Q_p, \mathbb V_{\hg \hh \mathbf{s}}) \longrightarrow H^1(\mathbb Q_p, \mathbb V_{\hg \hh \mathbf{s}}/\mathscr{F}^2\mathbb V_{\hg \hh \mathbf{s}}) \big).
\]
\end{definition}

Let $\mathrm{pr}^{-+}:H^1_{\bal}(\bb{Q}_p,\bb{V}_{\hg\hh \mathbf{s}})\rightarrow H^1(\bb{Q}_p,\bb{V}_{\hg}^-\hat{\otimes}\bb{V}_\hh^+\hat{\otimes}\Lambda(1-\mathbf{s}))$ be the map induced by the natural projection $\mathscr{F}^2\bb{V}_{\hg\hh\mathbf{s}}\rightarrow \bb{V}_\hg^-\hat{\otimes}\bb{V}_\hh^+\hat{\otimes}\Lambda(1-\mathbf{s})$ and let
\[
\mathcal{L}_{\hg\hh\mathbf{s}}^{\hg}:H^1_{\bal}(\bb{Q},\bb{V}_{\hg\hh\mathbf{s}})\longrightarrow I_{\hg}^{-1}I^{-1}\cl{O}_{\hg\hh\mathbf{s}}
\]
be the map defined by $\mathcal{L}_{\hg\hh\mathbf{s}}^\hg=\mathcal{L}_{\hg\hh\mathbf{s}}^{-+}\circ \mathrm{pr}^{-+}\circ \res_p$. Similarly we define $\mathcal{L}_{\hg\hh\mathbf{s}}^\hh$.

The construction of Beilinson--Flach classes was first carried out by Lei--Loeffler--Zerbes for fixed modular forms $(g,h)$ of weight two \cite{LLZ0}, and was later extended to Hida families $(\hg,\hh)$ in \cite{KLZ0}. Later, in \cite{LZ-coleman}, the construction was extended to the Coleman case; as discussed e.g. in \cite{LR1} there are situations where one must be more cautious and there may be some poles, for instance, at the critical $p$-stabilization of an Eisenstein series. However, in the setting considered in this note we will not deal with this kind of issues. 

Recall that, given a newform $\xi$ of level $N_{\xi}$, its image under the Atkin--Lehner operator $W_{N_\xi}$ is a scalar multiple of the conjugate eigenform $\xi^*$. Then, we define the {\it Atkin--Lehner pseudo-eigenvalue} $\lambda_{N_\xi}(\xi)$ of $\xi$ by $W_{N_\xi}(\xi) = \lambda_N(\xi)\xi^*$. As explained in \cite[\S10]{KLZ}, given a Hida family $\bm{\xi}$ of tame level $N_\xi$, there exists an element $\lambda_{N_\xi}(\bm{\xi})\in\cl{O}_{\bm{\xi}}^\times$ interpolating the Atkin--Lehner pseudo-eigenvalues of the crystalline specializations of $\bm{\xi}$.


We now state the main result of \cite{KLZ}.

\begin{theorem}\label{reclaw}
Fix an integer $c>1$ relatively prime to $6pN_gN_h$. Then, there exists a global cohomology class \[ _c\kappa(\hg,\hh,\mathbf{s}) \in H^1_{\bal}(\mathbb Q, \mathbb V_{\hg} \hat \otimes \mathbb V_{\hh} \hat \otimes \Lambda(1-\mathbf{s})) \] such that
\[
\cl{L}_{\hg\hh\mathbf{s}}^{\hg}(_c\kappa(\hg,\hh,\mathbf{s}))=\frac{(-1)^{\mathbf{s}}}{\lambda_{N_g}(\hg)} \cdot (c^2-c^{2\mathbf{s}-\mathbf{l}-\mathbf{m}+2}\chi_{\hg}^{(p)}(c)^{-1}\chi_{\hh}^{(p)}(c)^{-1}) \times L_p(\hg,\hh,\mathbf{s})
\]
and
\[
\cl{L}_{\hg\hh\mathbf{s}}^{\hh}(_c\kappa(\hg,\hh,\mathbf{s}))=\frac{(-1)^{\mathbf{s}}}{\lambda_{N_h}(\hh)} \cdot (c^2-c^{2\mathbf{s}-\mathbf{l}-\mathbf{m}+2}\chi_{\hg}^{(p)}(c)^{-1}\chi_{\hh}^{(p)}(c)^{-1}) \times L_p(\hh,\hg,\mathbf{s}).
\]
\end{theorem}
\begin{proof}
The global cohomology class $_c \kappa(\hg,\hh,\mathbf{s})$ is introduced in \cite[Definition 8.1.1]{KLZ}. The result follows from \cite[Proposition 8.1.7]{KLZ} and \cite[Theorem~B]{KLZ}.
\end{proof}


\begin{remark}\label{rk:dependenceofBFonc}
After tensoring with $\operatorname{Frac} (\mathcal{O}_{\mathbf{ghs}})$, the class
\[
\kappa(\hg, \hh,\mathbf{s}) \coloneqq C_c^{-1} \otimes {}_c \kappa(\hg, \hh,\mathbf{s})
\]
is independent of $c$, where
\begin{equation}
 \label{eq:defCd}
  C_c(\hg,\hh,{\bf s}) \coloneqq c^2 - c^{(2{\bf s}-{\bf l}-{\bf m}+2)}\chi_\hg^{(p)}(c)^{-1} \chi_\hh^{(p)}(c)^{-1}.
\end{equation}
\end{remark}


\subsection{Diagonal cycles}\label{subsec:diagonalcycles}

In this section, we recall the main results regarding the existence of $p$-adic families of diagonal cycles obtained in \cite{DR3} and \cite{BSV2} building on the previous geometric constructions of \cite{DR2}.

We use the notations and assumptions introduced in \S\ref{subsec:tripleproduct}. In particular, $(\hf,\hg,\hh)$ is a triple of Hida families of tame levels $(N_\hf,N_\hg,N_\hg)$ and characters $(\chi_\hf,\chi_\hg,\chi_\hh)$ and $(f,g,h)=(\hf_{x_0}^\circ,\hg_{y_0}^\circ,\hh_{z_0}^\circ)$ is a triple of good crystalline specializations of characters $(\chi_f,\chi_g,\chi_h)$ and weights $(k_0,l_0,m_0)$ with $k_{0}+l_{0}+m_{0}\equiv 0\pmod{2}$. 
We denote by $\mathbf{k}:\bb{Z}_p^\times \rightarrow \Lambda_\hf^\times$, $\mathbf{l}:\bb{Z}_p^\times\rightarrow \Lambda_\hg^\times$ and $\mathbf{m}:\bb{Z}_p^\times \rightarrow \Lambda_\hh^\times$ the corresponding tautological weight characters.

The running assumptions imply that $\chi_{\hf} \chi_{\hg} \chi_{\hh} = \omega^{2r}$ for some $r\in \mathbb{Z}$, where $\omega$ denotes the Teichm\"uller character. Thus we can choose a character ${\bf t}: \mathbb Z_p^\times \to \Lambda_{\hf\hg\hh}^\times=(\Lambda_\hf \hat{\otimes} \Lambda_\hg \hat{\otimes}\Lambda_\hh)^\times$ satisfying $2{\bf t} = {\bf k} + {\bf l} + {\bf m}$. There are two choices for $\mathbf{t}$, and we choose the one determined by imposing that the specialization of $\mathbf{t}$ at the fixed crystalline point $(x_0,y_0,z_0)$ is the character $\varepsilon_{\cyc}^{(k_0+l_0+m_0)/2}$. 
We also make the following assumption.

\begin{ass}
\hfill
    \begin{enumerate}[(i)]
        \item The Hida families $\hg$ and $\hh$ are residually irreducible and $p$-distinguished.
        \item $\bb{V}_{\hf}$ and $\bb{V}_{\hf}^-$ are free $\Lambda_\hf$-modules.
    \end{enumerate}
\end{ass}

Note that, by Proposition~\ref{prop:BSTW}, the second part of the assumption holds when $\hf$ is the CM Hida family introduced in \S\ref{subsec:wt1}, which is the case that we will consider in later sections.


We first introduce a Perrin-Riou logarithm that will be used for the formulation of an explicit reciprocity law. 

\begin{propo}\label{perrin2}
There exists a
homomorphism of $\Lambda_{\hf \hg \hh}$-modules
\[
\mathcal L_{\hf \hg \hh}^{+-+}: H^1(\mathbb Q_p, \bb{V}_\hf^+ \hat{\otimes}\mathbb V_{\hg}^- \hat \otimes \mathbb V_{\hh}^+ (2-\mathbf{t})) \rightarrow I_\hg^{-1} \cl{O}_{\hf \hg \hh },
\]
such that for every point $(x,y,z) \in \mathcal W_{{\bf fgh}}^{\mathrm{cl}}$ of weights $(k,l,m)$ with $\beta_{\hf_x}\alpha_{\hg_y}\beta_{\hh_z}\neq p^{(k+l+m-2)/2}$ the specialization of $\mathcal L_{\hf \hg \hh}^{+-+}$ at $(x,y,z)$ is the homomorphism
\[
\mathcal L_{{\bf fgh}}^{+-+}(x,y,z): H^1(\mathbb Q_p, V_{\hf_x}^+\otimes V_{\hg_y}^- \otimes V_{\hh_z}^+(1-c)) \rightarrow \mathbb C_p
\]
given by
\[
\mathcal L_{{\bf fgh}}^{+-+}(x,y,z) = \frac{1-p^{-c}\alpha_{\hf_x} \beta_{\hg_y}\alpha_{\hh_z}}{1-p^{-c} \beta_{\hf_x} \alpha_{\hg_y} \beta_{\hh_z}} \times \begin{cases} \frac{(-1)^{c-l}}{(c-l)!} \times \langle \log_{\BK}(\cdot), \omega_{\hf_x}\otimes\eta_{\hg_y}^\alpha \otimes \omega_{\hh_z} \rangle & \text{ if } l < k+ m \\ (l-c-1)! \times \langle \exp_{\BK}^*(\cdot), \omega_{\hf_x} \otimes\eta_{\hg_y}^\alpha \otimes \omega_{\hh_z} \rangle & \text{ if } l \geq k+m, \end{cases}
\]
where $\log_{\BK}$ denotes the Bloch-Kato logarithm, $\exp_{\BK}^*$ denotes the dual exponential map and $c=(k+l+m-2)/2$.
\end{propo}
\begin{proof}
This follows from \cite[Proposition~7.3]{BSV}.
\end{proof}


\begin{remark}
    We recall that the map in the previous theorem is obtained as the composition
    \[
    \mathcal L_{\hg \hh}^{+-+} = \langle \log_{\hf \hg \hh}^{+-+}(\cdot), \bm{\omega}_{\hf}\otimes\bm{\eta}_{\hg} \otimes \bm{\omega}_{\hh} \rangle,
    \]
    where $\log_{\hf \hg \hh}^{+-+}(\cdot)$ is the Perrin--Riou big logarithm introduced in \cite[Proposition~7.1]{BSV}. 
\end{remark}

\begin{remark}
    Letting $\hf$ (resp. $\hh$) play the role of $\hg$ in Proposition~\ref{perrin2}, we obtain a similar map $\cl{L}_{\hf\hg\hh}^{-++}$ (resp. $\cl{L}_{\hf\hg\hh}^{++-}$) with analogous properties.
\end{remark}

To shorten notation, put $\bb{V}_{\hf\hg\hh}^\dagger=\bb{V}_\hf\hat{\otimes}\bb{V}_\hg\hat{\otimes}\bb{V}_\hh(2-\mathbf{t})$. Let $\mathscr{F}^2\mathbb V_{\hf\hg\hh}^\dag\subset \mathbb V_{\hf\hg\hh}^{\dag}$ be the $G_{\bb{Q}_p}$-stable $\Lambda_{\hf \hg \hh}$-submodule of rank $4$ defined by
\[
\mathscr{F}^2\mathbb V_{\hf\hg\hh}^\dag=(\bb{V}_\hf\hat{\otimes}\bb{V}_{\hg}^+\hat{\otimes}\bb{V}_{\hh}^+ +\bb{V}_\hf^+\hat{\otimes}\bb{V}_{\hg}\hat{\otimes}\bb{V}_{\hh}^+ +\bb{V}_\hf^+\hat{\otimes}\bb{V}_{\hg}^+\hat{\otimes}\bb{V}_{\hh})(2-\mathbf{t}).
\]
Fix a finite set $\Sigma$ of places of $\mathbb Q$ containing $\infty$ and the primes dividing $N_fN_gN_hp$ and let $\mathbb Q^\Sigma$ be the maximal extension of $\mathbb Q$ unramified outside $\Sigma$.





\begin{definition}\label{def:Sel2}
The \textit{balanced Selmer group of $\bb{V}_{\hf \hg \hh}^\dagger$} is defined by
\[
H^1_{\bal}(\bb{Q},\mathbb V_{\hf \hg \hh}^{\dag}) = \ker \bigg(
H^1(\mathbb Q^\Sigma/\mathbb Q, \mathbb V_{\hf \hg \hh}^{\dag}) \longrightarrow
\frac{H^1(\mathbb Q_p,\mathbb V_{\hf \hg \hh}^\dag)}{H^1_{\bal}(\bb{Q}_p,\mathbb V_{\hf \hg \hh}^\dag)}\bigg),
\]
where
\[
H^1_{\bal}(\bb{Q}_p,\mathbb V_{\hf \hg \hh}^\dag)=\ker \big(H^1(\mathbb Q_p, \mathbb V_{\hf \hg \hh}^{\dag}) \longrightarrow H^1(\mathbb Q_p, \mathbb V_{\hf \hg \hh}^{\dag}/\mathscr{F}^2\mathbb V_{\hf \hg \hh}^{\dag}) \big).
\]
\end{definition}

Let $\mathrm{pr}^{+-+}:H^1_{\bal}(\bb{Q}_p,\bb{V}_{\hf\hg\hh}^\dagger)\rightarrow H^1(\bb{Q}_p,\bb{V}_{\hf}^+\hat{\otimes}\bb{V}_{\hg}^-\hat{\otimes}\bb{V}_\hh^+(2-\mathbf{t}))$ be the map induced by the natural projection $\mathscr{F}^2\bb{V}_{\hf\hg\hh}^\dagger\rightarrow \bb{V}_{\hf}^+\hat{\otimes}\bb{V}_\hg^-\hat{\otimes}\bb{V}_\hh^+(2-\mathbf{t})$ and let
\[
\mathcal{L}_{\hf\hg\hh}^{\hg}:H^1_{\bal}(\bb{Q},\bb{V}_{\hf\hg\hh}^\dagger)\longrightarrow I_{\hg}^{-1}\cl{O}_{\hf\hg\hh}
\]
be the map defined by $\mathcal{L}_{\hf\hg\hh}^\hg=\mathcal{L}_{\hf\hg\hh}^{+-+}\circ \mathrm{pr}^{+-+}\circ \res_p$. Similarly we define $\mathcal{L}_{\hf\hg\hh}^\hf$ and $\mathcal{L}_{\hf\hg\hh}^\hh$.

We now state main result of \cite{DR3} and \cite{BSV}.

\begin{theorem}\label{thm:reclawdiagonalcycles}
There exists a global cohomology class
\[
\kappa(\hf,\hg,\hh) \in H^1_{\bal}(\mathbb Q, \mathbb V_{\hf} \hat \otimes \mathbb V_{\hg} \hat \otimes \mathbb V_{\hh}(2-{\bf t}))
\]
such that, for $\bm{\xi}\in \lbrace \hf,\hg,\hh\rbrace$,
we have that
\[
\mathcal L_{\hf\hg\hh}^{\bm{\xi}}(\kappa(\hf, \hg, \hh)) = \Lp^{\xi}(\hf, \hg, \hh).
\]

\end{theorem}

\begin{proof}
This is \cite[Theorem~A]{BSV} or \cite[Theorem~5.1]{DR3}.

\end{proof}


\subsection{Anticyclotomic diagonal cycles}\label{subsec:ac-cycles}

We keep the notations and assumptions in the previous subsection, as well as the notations and assumptions in \S\ref{subsec:wt1}, and specialize the discussion in the previous subsection to the case in which $\hf$ is the Hida family introduced in \S\ref{subsec:wt1} with $f=\hf_{x_0}^\circ=\Eis_1(\varepsilon_K)$.

In this case, we have a class
\[
\kappa(\hf,\hg,\hh)\in H^1(\mathbb Q, \Ind_K^{\mathbb Q} \Lambda_{\hf}(\bm{\varphi}) \hat \otimes \mathbb V_{\hg} \hat \otimes \mathbb V_{\hh}(2-{\bf t}))
\]
which under Shapiro's isomorphism can also be seen as a class in
\[
H^1(K,\Lambda_\hf(\bm{\varphi})\hat{\otimes}\mathbb V_{\hg} \hat \otimes \mathbb V_{\hh}(2-\bf{t})).
\]
Moreover, after specializing $\hf$ to $f$, we obtain a class
\begin{align*}
    \kappa(f,\hg,\hh)\in H^1(\bb{Q}, (1\oplus\varepsilon_K)\otimes \bb{V}_\hg \hat{\otimes} \bb{V}_\hh (2-\mathbf{t}_1))\cong H^1(K,\bb{V}_\hg\hat{\otimes} \bb{V}_\hh(2-\bf{t}_1))
\end{align*}
where $\bf{t}_1$ is the composition of $\bf{t}$ with the map $(\Lambda_{\hf}\hat{\otimes}\Lambda_{\hg}\hat{\otimes}\Lambda_\hh)^\times\rightarrow(\Lambda_{\hg}\hat{\otimes}\Lambda_\hh)^\times$ induced by the point $x_0:\Lambda_\hf\rightarrow \cl{O}$.

Assume now in addition that $p$ does not divide the class number of $K$. Then the classes $\kappa(\hf,\hg,\hh)$ yield classes
\[
\tilde{\kappa}(\hf,\hg,\hh)\in H^1(K,\Lambda_{\ac}(\bm{\kappa}_{\ac}^{-1})\hat{\otimes} \bb{V}_\hg \hat{\otimes} \bb{V}_\hh (2-\mathbf{t}_1)) \cong H^1_{\mathrm{Iw}}(K_\infty^{\ac},\bb{V}_\hg\hat{\otimes}\bb{V}_{\hh}(2-\mathbf{t}_1)),
\]
where $K_\infty^{\ac}$ denotes the anticyclotomic $\bb{Z}_p$-extension of $K$, $\Lambda_{\ac}=\bb{Z}_p[[\Gal(K_\infty^{\ac}/K)]]$ denotes the corresponding Iwasawa algebra and $\bm{\kappa}_{\ac}:G_K\rightarrow \Lambda_{\ac}^\times$ denotes the tautological character. When $\hf$ is residually non-Eisenstein, \cite{ACR} realises the class $\tilde{\kappa}(\hf,\hg,\hh)$ as the bottom layer of an anticyclotomic Euler system. However, this does not apply to the choice of $\hf$ in this article. We intend to come back to this in the future.

\section{Comparison of classes}\label{sec:comparison}

We keep the assumptions and notations introduced in previous sections. In particular, $K$ is an imaginary quadratic field in which $p$ splits and $\hf$ is the CM Hida family introduced in \S\ref{subsec:wt1} with $f=\hf_{x_0}^\circ=\Eis_1(\varepsilon_K)$. Also, $(\hg,\hh)$ is a pair of Hida families of tame levels $(N_g,N_h)$ and characters $(\chi_\hg,\chi_\hh)$ such that $\chi_\hg \chi_\hh =\varepsilon_K \omega^{2r}$ for some $r\in \bb{Z}$, and $g=\hg_{y_0}^\circ$ and $h=\hh_{z_0}^\circ$ are good crystalline specializations of $\hg$ and $\hh$ of weights $l_0\geq 2$ and $m_0\geq 1$, respectively. They have characters $(\chi_g,\chi_h)=(\chi_{\hg}^{(p)},\chi_\hh^{(p)})$. We define $\alpha_g=a_p(\hg_{y_0})$, $\beta_g=\chi_g(p)p^{l_0-1}\alpha_g^{-1}$, $\alpha_h=a_p(\hh_{z_0})$, $\beta_h=\chi_h(p)p^{m_0-1}\alpha_h^{-1}$. The corresponding definitions for $f$ yield $\alpha_f=\beta_f=1$. Set $c_0=s_0=(l_0+m_0-1)/2$. The following set of assumptions will be in place in all this section.

\begin{ass}
\hfill
\begin{enumerate}[(i)]
    \item $\gcd(N_g,N_h)=1$;
    \item $\hg$ and $\hh$ are residually irreducible and $p$-distinguished;
    \item $h$ is a newform of level $N_h$ (i.e., $h$ has prime-to-$p$ conductor).
\end{enumerate}
\end{ass}


\subsection{An equality of $p$-adic $L$-functions}

Set $\Lambda_{\hg\hh}=\Lambda_{\hg}\hat{\otimes}\Lambda_{\hh}$ and $\cl{O}_{\hg \hh}=\Lambda_{\hg\hh}[1/p]$. Let $\Lp^g(f,\hg,\hh)\in I_{\hg}^{-1}\cl{O}_{\hg \hh}$ be the two-variable $p$-adic $L$-function obtained from the three-variable $p$-adic $L$-function $\Lp^g(\hf,\hg,\hh)$ introduced in \S\ref{subsec:tripleproduct} by specializing $\hf$ to $f$. Also, we define $L_p(\hg,\hh,(\mathbf{l}+\mathbf{m}-1)/2)\in I_\hg^{-1}\cl{O}_{\hg \hh}$ as the restriction of the three-variable $p$-adic $L$-function $L(\hg,\hh,\mathbf{s})$ introduced in \S\ref{subsec:Hida-Rankin} to the plane $2\mathbf{s}=\mathbf{l}+\mathbf{m}-1$. Similarly, we define the $p$-adic $L$-function $L_p(\hg,\hh\otimes\varepsilon_K,(\mathbf{l}+\mathbf{m}-1)/2)\in I_{\hg}^{-1}\cl{O}_{\hg \hh}$ by replacing $\hh$ by its twist $\hh\otimes \varepsilon_K$.

\begin{propo}\label{prop:factorization}
    We have the following equality of $p$-adic $L$-functions:
    \[
    \Lp^g(f,\hg,\hh)^2=L_p\left(\hg,\hh,\frac{\mathbf{l}+\mathbf{m}-1}{2}\right)L_p\left(\hg,\hh\otimes\varepsilon_K,\frac{\mathbf{l}+\mathbf{m}-1}{2}\right).
    \]
\end{propo}
\begin{proof}
    Since the points $(y,z)\in \tilde{\cl{W}}_{\hg}^{\mathrm{cl}}\times \tilde{\cl{W}}_{\hh}^{\mathrm{cl}}$ with weights $k_y>k_z$ and with $p\nmid \mathrm{cond}(\hg_y)\cdot \mathrm{cond}(\hh_z)$ form a Zariski dense subset of $\cl{U}_\hg\times\cl{U}_\hh$, it suffices to prove the equality after specialization at each such point. Note that for such a point we have that $(x_0,y,z)\in \cl{W}_{\hf \hg \hh}^g$.
    
    Fix such a point $(y,z)\in \tilde{\cl{W}}_\hg^{\mathrm{cl}}\times \tilde{\cl{W}}^{\mathrm{cl}}_\hh$ of weights $(l,m)$ with $l>m$ and let $c=(l+m-1)/2$. Since $V_f\cong L(1)\oplus L(\varepsilon_K)$, we deduce by Artin formalism the equality of complex $L$-values
    \[
    L(f,\hg_y^\circ,\hh_z^\circ,c)=L(\hg_y^\circ,\hh_z^\circ,c)L(\hg_y^\circ,\hh_z^\circ\otimes\varepsilon_K,c).
    \]
    Then, the equality
    \[
    \Lp^g(f,\hg,\hh)(y,z)^2=L_p\left(\hg,\hh,\frac{\mathbf{l}+\mathbf{m}-1}{2}\right)(y,z)\cdot L_p\left(\hg,\hh\otimes\varepsilon_K,\frac{\mathbf{l}+\mathbf{m}-1}{2}\right)(y,z)
    \]
    follows from Theorem~\ref{thm:Hida-3var} and Theorem~\ref{thm:hsieh}.
\end{proof}

From now on, we make the following assumption.

\begin{ass}\label{ass:nonvanishing}
    The $p$-adic $L$-function $\Lp^g(f,\hg,\hh)$ is not identically zero.
\end{ass}

\begin{remark}
    Note that, in light of Proposition~\ref{prop:factorization}, the assumption immediately implies that the $p$-adic $L$-functions $L_p(\hg,\hh,(\mathbf{l}+\mathbf{m}-1)/2)$ and $L_p(\hg,\hh\otimes\varepsilon_K,(\mathbf{l}+\mathbf{m}-1)/2)$ are not identically zero.
\end{remark}

\subsection{Weighted Beilinson--Flach class}\label{subsec:weightedBF}

As in previous sections, we denote by $\bb{V}_{\bm{\xi}}$ the big Galois representation associated with a Hida family $\bm{\xi}$. Note that the representations $\bb{V}_{\hh\otimes\varepsilon_K}$ and $\bb{V}_{\hh}\otimes \varepsilon_K$ are isomorphic. Since $p$ splits in $K$, an isomorphism of $\Lambda_{\hh}[G_{\bb{Q}}]$-modules $\bb{V}_{\hh\otimes\varepsilon_K}\cong \bb{V}_{\hh}\otimes \varepsilon_K$ yields an isomorphism of $\Lambda_{\hh}[G_{\bb{Q}_p}]$-modules $\bb{V}_{\hh\otimes\varepsilon_K}\cong \bb{V}_{\hh}$ and hence an isomorphism of Dieudonn\'e modules $\mathbf{D}(\bb{V}_{\hh\otimes\varepsilon_K}^+(1-\mathbf{m}-\chi_{\hh}))\cong \mathbf{D}(\bb{V}_{\hh}^+(1-\mathbf{m}-\chi_{\hh}))$. We choose an isomorphism $\iota:\bb{V}_{\hh\otimes\varepsilon_K}\rightarrow \bb{V}_\hh\otimes\varepsilon_K$ so that the composition
\[
\mathbf{D}(\bb{V}_{\hh\otimes\varepsilon_K}^+(1-\mathbf{m}-\chi_{\hh}))\xrightarrow{\iota} \mathbf{D}(\bb{V}_{\hh}^+(1-\mathbf{m}-\chi_{\hh})) \xrightarrow[]{\langle\, \cdot ,\,\bm{\omega}_\hh\rangle} \cl{O}_{\hh}
\]
agrees with the map
\[
\mathbf{D}(\bb{V}_{\hh\otimes\varepsilon_K}^+(1-\mathbf{m}-\chi_{\hh}))\xrightarrow[]{\langle\, \cdot ,\,\bm{\omega}_{\hh\otimes\varepsilon_K}\rangle} \cl{O}_{\hh}
\]
and from now on we identify $\bb{V}_{\hh\otimes\varepsilon_K}$ with $\bb{V}_\hh\otimes\varepsilon_K$ via this isomorphism.

We put $\bb{V}_{f\hg\hh}^\dagger=V_f\otimes\bb{V}_{\hg}\hat{\otimes} \bb{V}_\hh (2-\mathbf{t}_1)$ and $\bb{V}_{\hg\hh}^\dagger=\bb{V}_{\hg}\hat{\otimes}\bb{V}_\hh(2-\mathbf{t}_1)$. Note that
\[
\bb{V}_{f\hg\hh}^\dagger\cong(1\oplus \varepsilon_K)\otimes \bb{V}_{\hg\hh}^\dagger.
\]
After specializing to the plane $2\mathbf{s}=\mathbf{l}+\mathbf{m}-1$, the class $_c\kappa(\hg,\hh,\mathbf{s})$ introduced in \S\ref{subsec:BF} yields a class $_c\kappa(\hg,\hh)$ in $H^1(\bb{Q},\bb{V}_{\hg\hh}^\dagger)$. Replacing $\hh$ by $\hh\otimes\varepsilon_K$, we also obtain a class $_c\kappa(\hg,\hh\otimes\varepsilon_K)$ in $H^1(\bb{Q},\bb{V}_{\hg\hh}^\dagger\otimes\varepsilon_K)$. We define
\[
\kappa_{\hg,\hh}=\frac{(-1)^{\mathbf{s}}}{c^2-\chi_g^{-1}(c)\chi_h^{-1}(c)c} \cdot _c\kappa(\hg,\hh)\quad \text{and}\quad \kappa_{\hg,\hh\otimes\varepsilon_K}=\frac{(-1)^{\mathbf{s}}}{c^2-\chi_g^{-1}(c)\chi_h^{-1}(c)\varepsilon_K(c)c} \cdot _c\kappa(\hg,\hh\otimes\varepsilon_K).
\]
Note that, on account of Remark~\ref{rk:dependenceofBFonc}, the classes $\kappa_{\hg,\hh}$ and $\kappa_{\hg,\hh\otimes\varepsilon_K}$ do not depend on $c$. Note also that
\[
\kappa_{\hg,\hh}\in H^1_{\bal}(\bb{Q},\bb{V}_{\hg\hh}^\dagger) \quad \text{and}\quad \kappa_{\hg,\hh\otimes \varepsilon_K}\in H^1_{\bal}(\bb{Q},\bb{V}_{\hg\hh}^\dagger\otimes\varepsilon_K),
\]
where the Selmer groups $H^1_{\bal}(\bb{Q},\bb{V}_{\hg\hh}^\dagger)$ and $H^1_{\bal}(\bb{Q},\bb{V}_{\hg\hh}^\dagger\otimes\varepsilon_K)$ are defined as in \S\ref{subsec:BF}, specializing now all objects to the plane $2\mathbf{s}=\mathbf{l}+\mathbf{m}-1$.
Further, recall the basis $\{v_{f,1},v_{f,\varepsilon_K}\}$ of $V_f$ introduced in \S\ref{subsec:wt1}. We now use these elements to define the Beilinson--Flach class that we will use in the comparison.

\begin{defi}
The \emph{weighted Beilinson--Flach class} associated with the pair $(\hg,\hh)$ is the element
\[
\mathrm{BF}(f,\hg,\hh)=v_{f,1}\otimes L_p\left(\hg,\hh\otimes\varepsilon_K,\frac{\mathbf{l}+\mathbf{m}-1}{2}\right)\kappa_{\hg,\hh}+ v_{f,\varepsilon_K} \otimes L_p\left(\hg,\hh,\frac{\mathbf{l}+\mathbf{m}-1}{2}\right)\kappa_{\hg,\hh\otimes\varepsilon_K}
\]
in $H^1(K,\bb{V}_{f\hg\hh}^\dagger)$.   
\end{defi}
 
Since complex conjugation acts trivially on $\BF(f, \hg, \hh)$, according to the definitions of the vectors $v_{f,1}$ and $v_{f,\varepsilon_K}$ given in \ref{subsec:wt1}, 
this element descends to a class in $H^1(\mathbb Q, \bb{V}_{f\hg\hh}^\dagger)$ that we still denote with the same notation, i.e., we have
\[ 
\mathrm{BF}(f,\hg,\hh) \in H^1(\mathbb Q, \bb{V}_{f\hg \hh}^\dagger).
\]

\subsection{Selmer conditions}\label{subsec:Selmergroups}

Fix a finite set $\Sigma$ of places of $\mathbb Q$ containing $\infty$ and the primes dividing $N_fN_gN_hp$ and let $\mathbb Q^\Sigma$ be the maximal extension of $\mathbb Q$ unramified outside $\Sigma$.

Let $\mathscr{F}^2\bb{V}_{f\hg\hh}^\dagger\subset \bb{V}_{f\hg\hh}^\dagger$ be the $\Lambda_{\hg\hh}$-submodule of rank $4$ defined by
\[
\mathscr{F}^2\mathbb V_{f\hg\hh}^\dag=(\bb{V}_f\hat{\otimes}\bb{V}_{\hg}^+\hat{\otimes}\bb{V}_{\hh}^+ +\bb{V}_f^+\hat{\otimes}\bb{V}_{\hg}\hat{\otimes}\bb{V}_{\hh}^+ +\bb{V}_f^+\hat{\otimes}\bb{V}_{\hg}^+\hat{\otimes}\bb{V}_{\hh})(2-\mathbf{t}_1).
\]
Then, as in \S\ref{subsec:diagonalcycles}, we define the \textit{balanced Selmer group of $\bb{V}_{f\hg\hh}^\dagger$} by
\[
H^1_{\bal}(\bb{Q},\mathbb V_{f \hg \hh}^{\dag}) = \ker \bigg(
H^1(\mathbb Q^\Sigma/\mathbb Q, \mathbb V_{f \hg \hh}^{\dag}) \longrightarrow
\frac{H^1(\mathbb Q_p,\mathbb V_{f \hg \hh}^\dag)}{H^1_{\bal}(\bb{Q}_p,\mathbb V_{f \hg \hh}^\dag)}\bigg),
\]
where
\[
H^1_{\bal}(\bb{Q}_p,\mathbb V_{f \hg \hh}^\dag)=\ker \big(H^1(\mathbb Q_p, \mathbb V_{f \hg \hh}^{\dag}) \longrightarrow H^1(\mathbb Q_p, \mathbb V_{f \hg \hh}^{\dag}/\mathscr{F}^2\mathbb V_{f \hg \hh}^{\dag}) \big).
\]

We also need to introduce other Selmer conditions. Let $\bb{V}_{f\hg\hh}^g\subset \bb{V}_{f\hg\hh}^\dagger$ be the $\Lambda_{\hg\hh}$-submodule of rank $4$ defined by
\[
\bb{V}_{f\hg\hh}^g=V_f\otimes \bb{V}_{\hg}^+\hat{\otimes}\bb{V}_\hh(2-\mathbf{t}_1)
\]
and let $\bb{V}_{f\hg\hh}^{g\cup +}=\bb{V}_{f\hg\hh}^g+\mathscr{F}^2\bb{V}_{f\hg\hh}^\dagger$, which is a $\Lambda_{\hg\hh}$-submodule of $\bb{V}_{f\hg\hh}^\dagger$ of rank $5$.

\begin{defi}
For $\mathcal{L}\in\{\mathcal{G},\mathcal{G}\cup +\}$, the Selmer group $H^1_{\mathcal{L}}(\bb{Q},\mathbb V_{f\hg\hh}^{\dag})$ is defined by
\[
H^1_{\mathcal{L}}(\bb{Q},\mathbb V_{f\hg\hh}^{\dag}) = \ker \bigg(
H^1(\mathbb Q^\Sigma/\mathbb Q, \mathbb V_{f \hg\hh}^{\dag}) \longrightarrow
\frac{H^1(\mathbb Q_p,\mathbb V_{f \hg\hh}^\dag)}{H^1_{\mathcal{L}}(\mathbb Q_p,\mathbb V_{f \hg\hh}^\dag)}\bigg),
\]
where
\[
H_{\mathcal L}^1(\mathbb Q_p, \mathbb V_{f\hg\hh}^{\dag}) =
\begin{cases}
\ker \big(H^1(\mathbb Q_p, \mathbb V_{f\hg\hh}^{\dag}) \longrightarrow H^1(\mathbb Q_p, \mathbb V_{f \hg\hh}^{\dag}/\bb{V}_{f\hg\hh}^g) \big) &\textrm{if $\mathcal{L}=\cl{G}$},\\[0.5em]
\ker \big(H^1(\mathbb Q_p, \mathbb V_{f \hg\hh}^{\dag}) \longrightarrow H^1(\mathbb Q_p, \mathbb V_{f\hg\hh}^{\dag}/ \mathbb V_{f\hg\hh}^{g\cup +}) \big) &\textrm{if $\mathcal{L}=\mathcal{G}\cup +$}. 
\end{cases}
\]
\end{defi}

\subsection{Perrin-Riou maps}

In this subsection we introduce the Perrin-Riou maps that we will need for the comparison.

\begin{propo}\label{prop:2variablePR}
There exists a
homomorphism of $\Lambda_{\hg \hh}$-modules
\[
\mathcal L_{\hg \hh}^{-+}: H^1(\mathbb Q_p, \mathbb V_{\hg}^- \hat \otimes \mathbb V_{\hh}^+ (2-\mathbf{t}_1)) \rightarrow  I_{\hg}^{-1}\cl{O}_{\hg \hh},
\]
such that for every point $(y,z) \in \tilde{\cl{W}}_{\hg}^{\mathrm{cl}}\times \tilde{\cl{W}}_{\hh}^{\mathrm{cl}}$ of weights $(l,m)$ with $\alpha_{\hg_y}\beta_{\hh_z}\neq p^{(l+m-1)/2}$ the specialization of $\mathcal L_{\hg \hh}^{-+}$ at $(y,z)$ is the homomorphism
\[
\mathcal L_{{\bf gh}}^{-+}(y,z): H^1\left(\mathbb Q_p, V_{\hg_y}^- \otimes V_{\hh_z}^+\left(\frac{3-l-m}{2}\right)\right) \rightarrow \mathbb C_p
\]
given by
\[
\mathcal L_{{\bf gh}}^{-+}(y,z) = \frac{1-p^{(l+m-3)/2} \alpha_{\hg_y}^{-1} \beta_{\hh_z}^{-1}}{1-p^{(1-l-m)/2} \alpha_{\hg_y} \beta_{\hh_z}} \times \begin{cases} \frac{(-1)^{(m-l-1)/2}}{\left(\frac{m-l-1}{2}\right)!} \times \langle \log_{\BK}(\cdot), \eta_{\hg_y}^\alpha \otimes \omega_{\hh_z} \rangle & \text{ if } l < m+1 \\ \left(\frac{l-m-1}{2}\right)! \times \langle \exp_{\BK}^*(\cdot), \eta_{\hg_y}^\alpha \otimes \omega_{\hh_z} \rangle & \text{ if } l \geq m+1, \end{cases}
\]
where $\log_{\BK}$ is the Bloch--Kato logarithm and $\exp_{\BK}^*$ is the dual exponential map.
\end{propo}
\begin{proof}
    The existence and properties of the map $\cl{L}_{\hg\hh}^{-+}$ can be deduced from those of the map $\cl{L}_{\hg\hh\mathbf{s}}^{-+}$ in Proposition~\ref{perrin1} following the argument in \cite[Prop.~7.3]{BSV}.
\end{proof}

\begin{remark}
    Similarly, we have a $\Lambda_{\hg\hh}$-module homomorphism
    \[
    \mathcal L_{\hg \hh}^{+-}: H^1(\mathbb Q_p, \mathbb V_{\hg}^+ \hat \otimes \mathbb V_{\hh}^- (2-\mathbf{t}_1)) \rightarrow  I_{\hh}^{-1}\cl{O}_{\hg \hh},
    \]
    with analogous interpolation properties.
\end{remark}

\begin{lemma}
    The $\Lambda_{\hg\hh}$-module $H^1(\bb{Q}_p,\bb{V}_{\hg}^{-}\hat{\otimes}\bb{V}_{\hh}^+(2-\mathbf{t}_1))$ is torsion-free.
\end{lemma}
\begin{proof}
    To shorten notation, we write $\bb{V}=\bb{V}_{\hg}^{-}\hat{\otimes}\bb{V}_{\hh}^+(2-\mathbf{t}_1)$. Note that, as a $\Lambda_{\hg\hh}[G_{\bb{Q}_p}]$-module, we have that
    \[
    \bb{V}\cong \Lambda_{\hg\hh}(\underline{\chi}_g \underline{\alpha}_{\hg}^{-1} \underline{\alpha}_{\hh}\Theta),
    \]
    where $\underline{\chi}_g$ is the unramified character of $G_{\bb{Q}_p}$ defined by $\underline{\chi}_g(\Fr_p)=\chi_g(p)$, $\underline{\alpha}_{\hg}$ is the unramified character of $G_{\bb{Q}_p}$ defined by $\underline{\alpha}_{\hg}(\Fr_p)=a_p(\hg)$, $\underline{\alpha}_{\hh}$ is the unramified character of $G_{\bb{Q}_p}$ defined by $\underline{\alpha}_{\hh}(\Fr_p)=a_p(\hh)$ and $\Theta$ is the character of $G_{\bb{Q}_p}$ defined by
    \[
    \sigma\mapsto \omega_{\cyc}^{(l_0-m_0-1)/2}(\sigma) \langle \varepsilon_{\cyc}(\sigma)\rangle^{-1/2} [\langle \varepsilon_{\cyc}(\sigma)\rangle^{1/2}]_{\Lambda_{\hg}} \otimes [\langle \varepsilon_{\cyc}(\sigma)\rangle^{-1/2}]_{\Lambda_{\hh}}.
    \]
    Let $\lambda\in \Lambda_{\hg\hh}$. Then, we have an exact sequence
    \[
    H^0(\bb{Q}_p,\bb{V}/\lambda\bb{V})\rightarrow H^1(\bb{Q}_p,\bb{V})\rightarrow H^1(\bb{Q}_p,\bb{V}),
    \]
    where the second arrow is multiplication by $\lambda$. Thus, in order to show that the $\Lambda_{\hg\hh}$-module $H^1(\bb{Q}_p,\bb{V})$ is torsion-free, it suffices to show that $H^0(\bb{Q}_p,\bb{V}/\lambda\bb{V})=0$ for all non-zero $\lambda\in \Lambda_{\hg\hh}$. For that, it suffices to show that $H^0(\bb{Q}_p,\bb{V}/\frk{Q}\bb{V})=0$ for all height-$1$ prime ideal $\frk{Q}$ of $\Lambda_{\hg\hh}$. Note that, if $(l_0-m_0-1)/2\not\equiv 0 \pmod{p-1}$, then the inertia group $I_{\bb{Q}_p}$ acts non-trivially on $\bb{V}/\cl{I}\bb{V}$ for any proper ideal $\cl{I}$ of $\Lambda_{\hg\hh}$, so the result is immediate. Thus, we assume from now on that $(l_0-m_0-1)/2\equiv 0 \pmod{p-1}$ and therefore that $l_0\equiv m_0+1 \pmod{2(p-1)}$.

    Note that $\Lambda_{\hg\hh}$ is a finite flat extension of the unique factorization domain $\Lambda\hat{\otimes}\Lambda\simeq \bb{Z}_p[[X,Y]]$. Let $\gamma_0$ be a topological generator of $1+p\bb{Z}_p$ (e.g., we can take $\gamma_0=1+p$). If $\frk{q}$ is a height-one prime ideal of $\Lambda\hat{\otimes}\Lambda$ different from the ideal generated by $[\gamma_0]\otimes 1-1\otimes \gamma_0[\gamma_0]$, then the character $\Theta\vert_{I_{\bb{Q}_p}}:I_{\bb{Q}_p}\rightarrow ((\Lambda\hat{\otimes}\Lambda)/\frk{q})^\times$ is non-trivial and therefore $H^0(\bb{Q}_p,\bb{V}/\frk{Q}\bb{V})=0$ for any height-one prime ideal $\frk{Q}$ of $\Lambda_{\hg\hh}$ lying above $\frk{q}$.

    Now let $\frk{q}_1$ be the height-one prime ideal of $\Lambda\hat{\otimes}\Lambda$ generated by $[\gamma_0]\otimes 1-1\otimes \gamma_0[\gamma_0]$ and let $\frk{Q}_1$ be a height-one prime ideal of $\Lambda_{\hg\hh}$ above $\frk{q}_1$. Let $\frk{q}_2$ be the height-two prime ideal of $\Lambda\hat{\otimes}\Lambda$ corresponding to the arithmetic point $(P_{l-2,\bm{1}},P_{l-3,\bm{1}})\in\cl{W}\times\cl{W}$ for some integer $l>3$ such that $l\equiv l_0\pmod{2(p-1)}$. Note that $\frk{q}_1\subset \frk{q}_2$. Since $\Lambda_{\hg\hh}$ is a finite extension of $\Lambda\hat{\otimes}\Lambda$, we can find a prime ideal $\frk{Q}_2$ of $\Lambda_{\hg\hh}$ lying above $\frk{q}_2$ and such that $\frk{Q}_1\subset \frk{Q}_2$. The prime $\frk{Q}_2$ corresponds to a crystalline point $(y,z)\in \tilde{\cl{W}}_{\hg}^{\mathrm{cl}}\times \tilde{\cl{W}}_{\hh}^{\mathrm{cl}}$ of weights $(l,l-1)$. Since $l>3$, both $\hg_y$ and $\hh_z$ are $p$-old. Therefore, by the Ramanujan--Petersson conjecture, we have that the complex absolute value $\vert \alpha_{\hh_z}/\alpha_{\hg_y}\vert$ is $p^{-1/2}$. In particular $\chi_g(p)\alpha_{\hh_z}/\alpha_{\hg_y}\neq 1$. It follows that $H^0(\bb{Q}_p,\bb{V}/\frk{Q}_2\bb{V})=0$ and therefore $H^0(\bb{Q}_p,\bb{V}/\frk{Q}_1\bb{V})=0$.
\end{proof}

\begin{propo}\label{prop:injectivityofPR}
    The map $\cl{L}_{\hg\hh}^{-+}$ is injective.
\end{propo}
\begin{proof}
    By Theorem~\ref{reclaw},
    \[
    \cl{L}_{\hg\hh}^{-+}(\mathrm{pr}^{-+}(\res_p(\kappa_{\hg,\hh})))=\lambda_{N_g}(\hg)^{-1}\cdot L_p\left(\hg,\hh,\frac{\mathbf{l}+\mathbf{m}-1}{2}\right).
    \]
    In particular, it follows by Assumption~\ref{ass:nonvanishing} that the $\Lambda_{\hg\hh}$-module homomorphism $\cl{L}_{\hg\hh}^{-+}$ is non-zero. Since $H^1(\bb{Q}_p,\bb{V}_{\hg}^{-}\hat{\otimes}\bb{V}_{\hh}^+(2-\mathbf{t}_1))$ is a torsion-free $\Lambda_{\hg\hh}$-module of rank $1$, it immediately follows that $\cl{L}_{\hg\hh}^{-+}$ is injective.
\end{proof}

\begin{propo}\label{perrin2}
There exists a homomorphism of $\Lambda_{\hg \hh}$-modules
\[
\mathcal L_{f \hg \hh}^{+-+}: H^1(\mathbb Q_p, V_f^+ \otimes\mathbb V_{\hg}^- \hat \otimes \mathbb V_{\hh}^+ (2-\mathbf{t}_1)) \rightarrow I_\hg^{-1} \cl{O}_{\hg \hh },
\]
such that for every point $(y,z) \in \tilde{\cl{W}}_{\hg}^{\mathrm{cl}}\times \tilde{\cl{W}}_{\hh}^{\mathrm{cl}}$ of weights $(l,m)$ with $\alpha_{\hg_y}\beta_{\hh_z}\neq p^{(l+m-1)/2}$ the specialization of $\mathcal L_{f \hg \hh}^{+-+}$ at $(y,z)$ is the homomorphism
\[
\mathcal L_{{f\hg\hh}}^{+-+}(y,z): H^1\left(\mathbb Q_p, V_{f}^+\otimes V_{\hg_y}^- \otimes V_{\hh_z}^+\left(\frac{3-l-m}{2}\right)\right) \rightarrow \mathbb C_p
\]
given by
\[
\mathcal L_{{f\hg\hh}}^{+-+}(y,z) = \frac{1-p^{(l+m-3)/2} \alpha_{\hg_y}^{-1}\beta_{\hh_z}^{-1}}{1-p^{(1-l-m)/2} \alpha_{\hg_y} \beta_{\hh_z}} \times \begin{cases} \frac{(-1)^{(m-l-1)/2}}{\left(\frac{m-l-1}{2}\right)!} \times \langle \log_{\BK}(\cdot), \omega_{f}\otimes\eta_{\hg_y}^\alpha \otimes \omega_{\hh_z} \rangle & \text{ if } l < m+1 \\ \left(\frac{l-m-1}{2}\right)! \times \langle \exp_{\BK}^*(\cdot), \omega_{f} \otimes\eta_{\hg_y}^\alpha \otimes \omega_{\hh_z} \rangle & \text{ if } l \geq m+1, \end{cases}
\]
where $\log_{\BK}$ denotes the Bloch-Kato logarithm and $\exp_{\BK}^*$ denotes the dual exponential map.
\end{propo}
\begin{proof}
This follows as in \cite[Proposition~7.3]{BSV}, working with $V_f$ instead of $\bb{V}_{\hf}$.
\end{proof}

\begin{remark}
    Similarly, we can define maps $\cl{L}_{f\hg\hh}^{++-}$ and $\cl{L}_{f\hg\hh}^{-++}$.
\end{remark}

\begin{remark}
    The relation between $\cl{L}_{\hg\hh}^{-+}$ and $\cl{L}_{f\hg\hh}^{+-+}$ is as follows: given an element
    \[
    v\otimes z\in V_f^+\otimes H^1(\mathbb Q_p, \mathbb V_{\hg}^- \hat \otimes \mathbb V_{\hh}^+ (2-\mathbf{t}_1))\cong H^1(\mathbb Q_p, V_f^+ \otimes\mathbb V_{\hg}^- \hat \otimes \mathbb V_{\hh}^+ (2-\mathbf{t}_1)),
    \]
    we have that
    \[
    \cl{L}_{f\hg\hh}^{+-+}(v\otimes z)=\langle v,\omega_f\rangle \cdot \cl{L}_{\hg\hh}^{-+}(z).
    \]
\end{remark}

\subsection{The explicit comparison}\label{subsec:explicit_comparison}

We can now present one of the main results of this note, which was already anticipated in the introduction. The result gives a direct comparison between the weighted Beilinson--Flach class introduced in \S\ref{subsec:weightedBF} and a diagonal cycle class, in analogy with the main results of \cite{LR1}.

Recall that the class $\kappa_{\hg,\hh}$ belongs to the balanced Selmer group $H^1_{\bal}(\bb{Q},\bb{V}_{\hg\hh}^\dagger)$. Therefore, the image of $\res_p(\kappa_{\hg,\hh})$ in $H^1(\bb{Q}_p,\bb{V}_{\hg}^-\hat{\otimes}\bb{V}_{\hh}(2-\mathbf{t}_1))$ lands in $H^1(\bb{Q}_p,\bb{V}_{\hg}^-\hat{\otimes}\bb{V}_{\hh}^+(2-\mathbf{t}_1))$ and the image of $\res_p(\kappa_{\hg,\hh})$ in $H^1(\bb{Q}_p,\bb{V}_{\hg}\hat{\otimes}\bb{V}_{\hh}^-(2-\mathbf{t}_1))$ lands in $H^1(\bb{Q}_p,\bb{V}_{\hg}^+\hat{\otimes}\bb{V}_{\hh}^-(2-\mathbf{t}_1))$. We denote these images by $\res_p(\kappa_{\hg,\hh})^{-+}$ and $\res_p(\kappa_{\hg,\hh})^{+-}$, respectively. Similarly we define $\res_p(\kappa_{\hg,\hh\otimes\varepsilon_K})^{-+}$ and $\res_p(\kappa_{\hg,\hh\otimes\varepsilon_K})^{+-}$.

We now establish the key proposition towards the comparison that we will prove later.

\begin{proposition}\label{prop:localconditionforBF}
    The image of $\res_p(\BF(f, \hg, \hh))$ in $H^1(\bb{Q}_p, V_f \otimes \bb{V}_{\hg}^- \hat{\otimes} \bb{V}_{\hh}(2-\mathbf{t}_1))$ lands in $H^1(\bb{Q}_p,V_f^+\otimes \bb{V}_\hg^- \hat{\otimes} \bb{V}_\hh^+(2-\mathbf{t}_1))$ and is given by
    \[
    v_f^+\otimes \left(L_p\left(\hg,\hh\otimes\varepsilon_K,\frac{\mathbf{l}+\mathbf{m}-1}{2}\right)\res_p(\kappa_{\hg,\hh})^{-+}+L_p\left(\hg,\hh,\frac{\mathbf{l}+\mathbf{m}-1}{2}\right)\res_p(\kappa_{\hg,\hh\otimes\varepsilon_K})^{-+}\right)
    \]
\end{proposition}
\begin{proof}
    Since $\kappa_{\hg,\hh}$ belongs to $H^1_{\bal}(\bb{Q},\bb{V}_{\hg\hh}^\dagger)$ and $\kappa_{\hg,\hh\otimes\varepsilon_K}$ belongs to $H^1_{\bal}(\bb{Q},\bb{V}_{\hg\hh}^\dagger\otimes\varepsilon_K)$, the image of $\res_p(\BF(f, \hg, \hh))$ in $H^1(\bb{Q}_p, V_f \otimes \bb{V}_{\hg}^- \hat{\otimes} \bb{V}_{\hh}(2-\mathbf{t_1}))$ lies in $H^1(\bb{Q}_p, V_f \otimes \bb{V}_{\hg}^- \hat{\otimes} \bb{V}_{\hh}^+(2-\mathbf{t_1}))$. Therefore, to conclude the proof it suffices to show that the class
    \[
    v_f^-\otimes \left(L_p\left(\hg,\hh\otimes\varepsilon_K,\frac{\mathbf{l}+\mathbf{m}-1}{2}\right)\res_p(\kappa_{\hg,\hh})^{-+}-L_p\left(\hg,\hh,\frac{\mathbf{l}+\mathbf{m}-1}{2}\right)\res_p(\kappa_{\hg,\hh\otimes\varepsilon_K})^{-+}\right)
    \]
    in $H^1(\bb{Q}_p, V_f^- \otimes \bb{V}_{\hg}^- \hat{\otimes} \bb{V}_{\hh}^+(2-\mathbf{t_1}))$ is zero. Using Theorem~\ref{reclaw}, we have that
    \begin{align*}
        \mathcal{L}_{\hg \hh}^{-+} &\left(L_p\left(\hg,\hh\otimes\varepsilon_K,\frac{\mathbf{l}+\mathbf{m}-1}{2}\right)\res_p(\kappa_{\hg,\hh})^{-+}-L_p\left(\hg,\hh,\frac{\mathbf{l}+\mathbf{m}-1}{2}\right)\res_p(\kappa_{\hg,\hh\otimes\varepsilon_K})^{-+}\right) \\
        &= \lambda_{N_g}(\hg)^{-1}\cdot L_p\left(\hg,\hh\otimes\varepsilon_K,\frac{\mathbf{l}+\mathbf{m}-1}{2}\right)L_p\left(\hg,\hh,\frac{\mathbf{l}+\mathbf{m}-1}{2}\right) \\ &-\lambda_{N_g}(\hg)^{-1}\cdot L_p\left(\hg,\hh,\frac{\mathbf{l}+\mathbf{m}-1}{2}\right)L_p\left(\hg,\hh\otimes\varepsilon_K,\frac{\mathbf{l}+\mathbf{m}-1}{2}\right)=0.
    \end{align*}
    Since the map $\cl{L}_{\hg \hh}^{-+}$ is injective by Proposition~\ref{prop:injectivityofPR}, the result follows.
\end{proof}

\begin{remark}\label{rk:localconditionforBF}
    As a consequence of the previous proposition, we have that
    \[
    \BF(f,\hg,\hh)\in H^1_{\cl{G}\cup+}(\bb{Q},\bb{V}_{f\hg\hh}^\dagger).
    \]
\end{remark}

Let $\kappa(f,\hg,\hh)\in H^1_{\bal}(\bb{Q},\bb{V}_{f\hg\hh}^\dagger)$ be the class obtained from the class $\kappa(\hf,\hg,\hh)\in H^1_{\bal}(\bb{Q},\bb{V}_{\hf\hg\hh}^\dagger)$ introduced in Theorem~\ref{thm:reclawdiagonalcycles} by specializing $\hf$ to $f$. We denote by $\res_p(\kappa(f,\hg,\hh))^{+-+}$ the image of the class $\res_p(\kappa(f,\hg,\hh))$ in $H^1(\bb{Q}_p,V_f^+\otimes\bb{V}_\hg^-\hat{\otimes}\bb{V}_\hh^+(2-\mathbf{t}_1))$. Similarly, we denote by $\res_p(\BF(f,\hg,\hh))^{+-+}$ the image of $\res_p(\BF(f,\hg,\hh))$ in $H^1(\bb{Q}_p,V_f^+\otimes\bb{V}_\hg^-\hat{\otimes}\bb{V}_\hh^+(2-\mathbf{t}_1))$.

We also introduce the following definition.

\begin{defi}\label{def:period}
Let $\Omega_{f, \gamma} \in L^{\times}$ be the $p$-adic period given by
\[
\Omega_{f,\gamma} = 2 \cdot \langle v_f^+, \omega_f \rangle_f.
\] 
\end{defi}

\begin{remark}
The definition of $\Omega_{f,\gamma}$ depends on the isomorphism $\gamma$ fixed in \S\ref{subsec:wt1}, hence the notation.    
\end{remark}

\begin{theorem}\label{thm:main}
    The equality
    \[
    \Omega_{f,\gamma} \cdot \cl{L}_{f\hg\hh}^{+-+}(\res_p(\kappa(f,\hg,\hh))^{+-+})^2 =\lambda_{N_g}(\hg) \cdot \mathcal{L}_{f\hg\hh}^{+-+}(\res_p(\BF(f, \hg, \hh))^{+-+})
    \]
    holds.
\end{theorem}
\begin{proof}
    By Proposition~\ref{prop:localconditionforBF}, we have that
    \begin{align*}
        \mathcal{L}^{+-+}_{f\hg\hh}(\res_p(\BF(f,\hg,\hh))^{+-+})&=\langle v_f^+,\omega_f\rangle \cdot L_p\left(\hg,\hh\otimes\varepsilon_K,\frac{\mathbf{l}+\mathbf{m}-1}{2}\right)\cl{L}_{\hg\hh}^{-+}\left(\res_p(\kappa_{\hg,\hh})^{-+}\right) \\
        &+\langle v_f^+,\omega_f\rangle \cdot L_p\left(\hg,\hh,\frac{\mathbf{l}+\mathbf{m}-1}{2}\right)\cl{L}_{\hg\hh}^{-+}\left(\res_p(\kappa_{\hg,\hh\otimes\varepsilon_K})^{-+}\right).
    \end{align*}
    Using Theorem~\ref{reclaw}, it follows that
    \[
    \mathcal{L}^{+-+}_{f\hg\hh}(\res_p(\BF(f,\hg,\hh))^{+-+})=\frac{\Omega_{f,\gamma}}{\lambda_{N_g}(\hg)}\cdot L_p\left(\hg,\hh,\frac{\mathbf{l}+\mathbf{m}-1}{2}\right)L_p\left(\hg,\hh\otimes\varepsilon_K,\frac{\mathbf{l}+\mathbf{m}-1}{2}\right).
    \]
    Hence, by Proposition~\ref{prop:factorization},
    \[
    \lambda_{N_g}(\hg)\cdot \mathcal{L}^{+-+}_{f\hg\hh}(\res_p(\BF(f,\hg,\hh))^{+-+})=\Omega_{f,\gamma}\cdot \Lp^g(f,\hg,\hh)^2,
    \]
    which is equal to
    \[
    \Omega_{f,\gamma} \cdot \cl{L}^{+-+}_{f\hg\hh}(\res_p(\kappa(f,\hg,\hh))^{+-+})^2
    \]
    by Theorem~\ref{thm:reclawdiagonalcycles}.
\end{proof}


\begin{corollary}\label{corollary:comparison}
    Assume that $H^1_{\cl{G}\cup +}(\bb{Q}, \bb{V}_{f\hg\hh}^\dagger)$ is a torsion-free $\Lambda_{\hg\hh}$-module of rank $1$. Then $\BF(f,\hg,\hh)$ belongs to $H^1_{\bal}(\bb{Q}, \bb{V}_{f\hg\hh}^\dagger)$ and
    \[
    \lambda_{N_g}(\hg)\cdot\BF(f,\hg,\hh)=\Omega_{f,\gamma}\cdot \Lp^g(f,\hg,\hh) \cdot \kappa(f,\hg,\hh).
    \]

\end{corollary}
\begin{proof}
    Let $\mathrm{pr}^{+-+}:H^1_{\cl{G}\cup +}(\bb{Q}_p,\bb{V}_{f\hg\hh}^\dagger)\rightarrow H^1(\bb{Q}_p,V_f^+\otimes\bb{V}_\hg^-\hat{\otimes}\bb{V}_\hh^+(2-\mathbf{t}_1))$ be the map induced by the natural projection $\bb{V}_{f\hg\hh}^{g\cup +}\rightarrow V_f^+\otimes\bb{V}_\hg^-\hat{\otimes}\bb{V}_\hh^+(2-\mathbf{t}_1)$. Since the element $\Lp^g(f,\hg,\hh)\in I_{\hg}^{-1}\cl{O}_{\hg\hh}$ is non-zero by Assumption~\ref{ass:nonvanishing}, it follows from Theorem~\ref{thm:reclawdiagonalcycles} that the map
    \[
    \cl{L}_{f\hg\hh}^{+-+}\circ\mathrm{pr}^{+-+}\circ\res_p: H^1_{\cl{G}\cup +}(\bb{Q},\bb{V}_{f\hg\hh}^\dagger) \longrightarrow I_{\hg}^{-1} \cl{O}_{\hg\hh}
    \]
    is non-zero and therefore injective by our assumptions on $H^1_{\cl{G}\cup +}(\bb{Q}, \bb{V}_{f\hg\hh}^\dagger)$. Since both $\kappa(f,\hg,\hh)$ and $\BF(f,\hg,\hh)$ belong to $H^1_{\cl{G}\cup +}(\bb{Q},\bb{V}_{f\hg\hh}^\dagger)$, the result now follows immediately from the previous theorem. 
\end{proof}

\begin{remark}\label{rk:torsionfreeness}
    Regarding the assumptions in the previous corollary, it is expected by sign considerations that the $\Lambda_{\hg\hh}$-module $H^1_{\bal}(\bb{Q}, \bb{V}_{f\hg\hh}^\dagger)$ has rank $1$ and that the $\Lambda_{\hg\hh}$-module $H^1_{\cl{G}}(\bb{Q}, \bb{V}_{f\hg\hh}^\dagger)$ has rank zero, which would easily imply that the $\Lambda_{\hg\hh}$-module $H^1_{\cl{G}\cup +}(\bb{Q},\bb{V}_{f\hg\hh}^\dagger)$ has rank $1$. Moreover, we can ensure that $H^1_{\cl{G}\cup +}(\bb{Q}, \bb{V}_{f\hg\hh}^\dagger)$ is torsion-free by imposing the condition that $H^0(\bb{Q},\overline{\rho}^\dagger)=0$, where $\overline{\rho}^\dagger$ is the residual representation attached to $\bb{V}_{f\hg \hh}^\dagger$.
\end{remark}

\begin{remark}
It may be instructive to discuss the analogies and differences between a comparison of this kind and that carried out in \cite{LR1}, where we also used a Coleman family passing through a critical Eisenstein series. Here, the main idea is to construct cohomology classes over $K$, that we can then descend and compare with a class over $\mathbb Q$. However, in \cite{LR1}, the classes are all defined over $\mathbb Q$, so the comparison is between a suitable projection of one of the classes and the other.
Moreover, weight-one modular forms behave in an ostensibly different way, since both $p$-stabilizations are ordinary. 
\end{remark}

\begin{corollary}
    Assume that $H^1_{\cl{G}\cup +}(\bb{Q}, \bb{V}_{f\hg\hh}^\dagger)$ is a torsion-free $\Lambda_{\hg\hh}$-module of rank $1$. Then
    \begin{align*}
     \Lp^g(f,\hg,\hh)\Lp^h(f,\hg,\hh)&=\frac{\lambda_{N_g}(\hg)}{\lambda_{N_h}(\hh)}\cdot L_p\left(\hg,\hh,\frac{\mathbf{l}+\mathbf{m}-1}{2}\right)L_p\left(\hh\otimes\varepsilon_K,\hg,\frac{\mathbf{l}+\mathbf{m}-1}{2}\right) \\ &=\frac{\lambda_{N_g}(\hg)}{\lambda_{N_h}(\hh)}\cdot L_p\left(\hg,\hh\otimes\varepsilon_K,\frac{\mathbf{l}+\mathbf{m}-1}{2}\right)L_p\left(\hh,\hg,\frac{\mathbf{l}+\mathbf{m}-1}{2}\right).
    \end{align*}
\end{corollary}
\begin{proof}
    By the previous corollary, we know that $\BF(f,\hg,\hh)\in H^1_{\bal}(\bb{Q},\bb{V}_{f\hg\hh}^\dagger)$. Therefore, the image of $\BF(f,\hg,\hh)$ in $H^1(\bb{Q}_p,V_f^-\otimes\bb{V}_{\hg}\hat{\otimes}\bb{V}_{\hh}^{-}(2-\mathbf{t}_1))$ is zero. This image is given by
    \[
    v_f^-\otimes \left(L_p\left(\hg,\hh\otimes\varepsilon_K,\frac{\mathbf{l}+\mathbf{m}-1}{2}\right)\res_p(\kappa_{\hg,\hh})^{+-}-L_p\left(\hg,\hh,\frac{\mathbf{l}+\mathbf{m}-1}{2}\right)\res_p(\kappa_{\hg,\hh\otimes\varepsilon_K})^{+-}\right),
    \]
    so it follows that
    \[
    L_p\left(\hg,\hh\otimes\varepsilon_K,\frac{\mathbf{l}+\mathbf{m}-1}{2}\right)\res_p(\kappa_{\hg,\hh})^{+-}-L_p\left(\hg,\hh,\frac{\mathbf{l}+\mathbf{m}-1}{2}\right)\res_p(\kappa_{\hg,\hh\otimes\varepsilon_K})^{+-}=0.
    \]
    Applying the Perrin-Riou map $\mathcal{L}_{\hg\hh}^{+-}$, we deduce by Theorem~\ref{reclaw} that
    \begin{align}
    &L_p\left(\hg,\hh,\frac{\mathbf{l}+\mathbf{m}-1}{2}\right)L_p\left(\hh\otimes\varepsilon_K,\hg,\frac{\mathbf{l}+\mathbf{m}-1}{2}\right) \nonumber \\
    &=L_p\left(\hg,\hh\otimes\varepsilon_K,\frac{\mathbf{l}+\mathbf{m}-1}{2}\right)L_p\left(\hh,\hg,\frac{\mathbf{l}+\mathbf{m}-1}{2}\right). \label{eq:equalityofHRpadicLfunctions}
    \end{align}
    Now note that, for any element
    \[
    v\otimes z \in V_f^+ \otimes H^1(\bb{Q}_p,\bb{V}_\hg^+\hat{\otimes} \bb{V}_\hh^-(2-\mathbf{t}_1))\cong H^1(\bb{Q}_p,V_f^+\otimes \bb{V}_{\hg}^+\hat{\otimes}\bb{V}_{\hh}^-(2-\mathbf{t}_1)),
    \]
    we have that
    \[
    \mathcal{L}_{f\hg\hh}^{++-}(v\otimes z)=\langle v,\omega_f\rangle \cdot\mathcal{L}_{\hg\hh}^{+-}(z).
    \]
    Let $\mathrm{pr}^{++-}:H^1_{\bal}(\bb{Q}_p,\bb{V}_{f\hg\hh}^\dagger)\rightarrow H^1(\bb{Q}_p,V_f^+\otimes\bb{V}_\hg^+\hat{\otimes}\bb{V}_\hh^-(2-\mathbf{t}_1))$ be the map induced by the natural projection $\mathscr{F}^2\bb{V}_{f\hg\hh}^\dagger\rightarrow V_f^+\otimes\bb{V}_\hg^+\hat{\otimes}\bb{V}_\hh^-(2-\mathbf{t}_1)$. Using Theorem~\ref{reclaw} and equation~\ref{eq:equalityofHRpadicLfunctions}, the image of $\BF(f,\hg,\hh)$ by $\cl{L}_{f\hg\hh}^{++-}\circ \mathrm{pr}^{++-}\circ \res_p$ is equal to
    \[
    \Omega_{f,\gamma}\cdot\lambda_{N_h}(\hh)^{-1}\cdot L_p\left(\hg,\hh,\frac{\mathbf{l}+\mathbf{m}-1}{2}\right)L_p\left(\hh\otimes\varepsilon_K,\hg,\frac{\mathbf{l}+\mathbf{m}-1}{2}\right).
    \]
    Also, by Theorem~\ref{thm:reclawdiagonalcycles}, the image of $\kappa(f,\hg,\hh)$ by $\cl{L}_{f\hg\hh}^{++-}\circ \mathrm{pr}^{++-}\circ \res_p$ is equal to $\Lp^h(f,\hg,\hh)$. Hence, the result now follows from the previous corollary.
    
\end{proof}

\begin{remark}
    Note that the $p$-adic $L$-functions involved in the statement have disjoint interpolation ranges, so the result does not follow by a direct comparison of complex $L$-values.
\end{remark}

\subsection{Formulae for specializations}\label{subsec:nonexceptional}

We can now specialize the Hida families $(\hg,\hh)$ and obtain similar results for the fixed modular forms $(g,h)$. We write $(g_\alpha,h_\alpha)=(\hg_{y_0},\hh_{z_0})$.

Let $\cl{L}_{g_\alpha h_\alpha}^{-+}:H^1(\bb{Q}_p,V_g^-\otimes V_h^+(1-c_0))\rightarrow \bb{C}_p$ and $\cl{L}_{fg_\alpha h_\alpha}^{+-+}:H^1(\bb{Q}_p,V_f^+\otimes V_g^-\otimes V_h^+(1-c_0))\rightarrow \bb{C}_p$ be the maps obtained from the maps $\cl{L}_{\hg\hh}^{+-}$ and $\cl{L}_{f\hg\hh}^{+-+}$ introduced above by specializing $(\hg,\hh)$ to $(g_\alpha,h_\alpha)$. Note that, for any element
\[
v\otimes z \in V_f^+ \otimes H^1(\bb{Q}_p,V_g^-\otimes V_h^+(1-c_0))\cong H^1(\bb{Q}_p,V_f^+\otimes V_g^-\otimes V_h^+(1-c_0)),
\]
we have that
\[
\mathcal{L}_{fg_\alpha h_\alpha}^{+-+}(v\otimes z)=\langle v,\omega_f\rangle \cdot\mathcal{L}_{gh}^{-+}(z).
\]
Also note that, if $\alpha_g \beta_h\neq p^{l_0+m_0-1}$, then we have
\[
\cl{L}_{g_\alpha h_\alpha}^{-+}= \frac{1-p^{(l_0+m_0-3)/2} \alpha_{g}^{-1} \beta_{h}^{-1}}{1-p^{(-l_0-m_0+1)/2} \alpha_{g} \beta_{h}} \times \begin{cases} \frac{(-1)^{(m_0-l_0-1)/2}}{\left(\frac{m_0-l_0-1}{2}\right)!} \times \langle \log_{\BK}(\cdot), \eta_{g_\alpha}^\alpha \otimes \omega_{h_\alpha} \rangle & \text{ if } l_0 \leq m_0 \\ \left(\frac{l_0-m_0-1}{2}\right)! \times \langle \exp_{\BK}^*(\cdot), \eta_{g_\alpha}^\alpha \otimes \omega_{h_\alpha} \rangle & \text{ if } l_0 > m_0. \end{cases}
\]

To shorten notation, we write $V_{fgh}^\dagger=V_f\otimes V_g\otimes V_h(1-c_0)$ and we define Selmer conditions and Selmer groups for this representation analogous to the ones introduced in \S\ref{subsec:Selmergroups} for the representation $\bb{V}_{f\hg\hh}^\dagger$. We denote by $\kappa(f,g_\alpha,h_\alpha)\in H^1_{\bal}(\bb{Q},V_{fgh}^\dagger)$ the specialization of the class $\kappa(f,\hg,\hh)$ and we denote by $\BF(f,g_\alpha,h_\alpha)\in H^1_{\cl{G}\cup +}(\bb{Q},V_{fgh}^\dagger)$ the specialization of the class $\BF(f,\hg,\hh)$.

\begin{proposition}\label{prop:main-nonexceptional}
    The equality
    \[
    \Omega_{f,\gamma} \cdot \cl{L}_{fg_\alpha h_\alpha}^{+-+}(\res_p(\kappa(f,g_\alpha,h_\alpha))^{+-+})^2 =\lambda_{N_g}(g) \cdot \mathcal{L}_{f g_\alpha h_\alpha}^{+-+}(\res_p(\BF(f, g_\alpha, h_\alpha))^{+-+})
    \]
    holds.
\end{proposition}
\begin{proof}
    This is an immediate consequence of Theorem~\ref{thm:main}.
\end{proof}

\begin{corollary}\label{cor:non-exceptionalcomparison}
    Assume that $H^1_{\bal}(\bb{Q}, V_{fgh}^\dagger)$ is $1$-dimensional and that $\Lp^g(f,g_\alpha,h_\alpha)\neq 0$. If $l_0=2$ and $m_0=1$, further assume that $\alpha_g\beta_h\neq 1$. Then $\BF(f,g_\alpha,h_\alpha)$ belongs to $H^1_{\bal}(\bb{Q}, V_{fgh}^\dagger)$ and
    \[
    \lambda_{N_g}(g)\cdot\BF(f,g_\alpha,h_\alpha)=\Omega_{f,\gamma}\cdot \Lp^g(f,g_\alpha,h_\alpha) \cdot \kappa(f,g_\alpha,h_\alpha).
    \]
\end{corollary}
\begin{proof}
    Arguing as in the proof of \cite[Thm~9.5]{ACR}, we have $H^1_{\cl{G}\cup+}(\bb{Q},V_{fgh}^\dagger)=H^1_{\bal}(\bb{Q},V_{fgh}^\dagger)$. Now the proof follows as in Corollary~\ref{corollary:comparison}.
\end{proof}

\begin{remark}
We emphasize that the class $\kappa(f,g_\alpha,h_\alpha)$ can be interpreted as an anticyclotomic cohomology class via the identifications discussed in \S\ref{subsec:ac-cycles}.
\end{remark}

\subsection{A factorization formula for the big logarithm of a Beilinson--Flach class}

As a consequence of the cyclotomic results of \cite{BDV}, B\"uy\"ukboduk, Casazza, and Sakamoto \cite{BC}, \cite{BS} obtained an expression of the Hida--Rankin $p$-adic $L$-function in terms of the Ochiai big logarithm of the Kato class. Here, we discuss an analogue of \cite[Prop. 8.9]{BC}, showing that the image of $\BF(f,\hg,\hh)$ under a suitable Perrin-Riou map factors as the product of two triple product $p$-adic $L$-functions. 

Let $\cl{L}_{f\hg\hh}^{-++}:H^1(\bb{Q}_p,V_f^-\otimes \bb{V}_{\hg}^+\hat{\otimes}\bb{V}_{\hh}^+(2-\mathbf{t}_1))\rightarrow \cl{O}_{\hg\hh}$ be the Perrin-Riou map obtained from $\cl{L}_{\hf\hg\hh}$ by specializing $\hf$ to $f$. It is given by the composition of a big logarithm map $\log_{f \hg \hh}^{-++}$ with the map obtained by pairing with $\eta_{f_\alpha}^\alpha\otimes \bm{\omega}_{\hg}\otimes \bm{\omega}_{\hh}$. Note that the pairing is taken with respect to the differential $\eta_{f_\alpha}^\alpha$ corresponding to the $p$-stabilized form $f_\alpha$, as opposed to $\eta_f^\alpha$.

\begin{defi}\label{def:big-log}
We denote by $\Log_{\bm{\omega}_{\hg} \otimes \bm{\omega}_{\hh}}$ the $\Lambda_{\hg\hh}$-module homomorphism
\[
\Log_{\bm{\omega}_{\hg} \otimes \bm{\omega}_{\hh}}=\cl{L}^{-++}_{f\hg\hh}\circ\mathrm{pr}^{-++}\circ\res_p : H^1_{\bal}(\bb{Q},\bb{V}_{f\hg\hh}^\dagger)\longrightarrow  \cl{O}_{\hg\hh}.
\]
\end{defi}

\begin{remark}\label{rk:other-proj}
In terms of the Beilinson--Flach classes, we can think that the pairing with $\eta_{\hg} \otimes \omega_{\hh}$ (resp. $\omega_{\hg} \otimes \eta_{\hh}$) allows us to recover the $p$-adic $L$-function $L_p(\hg,\hh)$ (resp. $L_p(\hh,\hg)$), while here we are considering instead the pairing with $\omega_{\hg} \otimes \omega_{\hh}$.
\end{remark}

Let $\Lp^f(f,\hg,\hh)\in \cl{O}_{\hg\hh}$ be the two-variable $p$-adic $L$-function obtained from the three-variable $f$-dominant $p$-adic $L$-function $\Lp^f(\hf,\hg,\hh)$ introduced in \S\ref{subsec:tripleproduct} by specializing $\hf$ to $f$. In this subsection, we relate the $p$-adic $L$-function $\Lp^f(f,\hg,\hh)$ to the weighted Beilinson--Flach class $\BF(f,\hg,\hh)$. 
Note that, while $\Lp^g(f, \hg, \hh)^2$ and $\Lp^h(f, \hg, \hh)^2$ may be factored as the product of two Hida--Rankin $p$-adic $L$-functions by $p$-adic Artin formalism, 
the case where the dominant $p$-adic $L$-function is Eisenstein is subtler, as already shown in \cite{BC}. 

\begin{proposition}
Assume that $H^1_{\cl{G}\cup +}(\bb{Q}, \bb{V}_{f\hg\hh}^\dagger)$ is a torsion-free $\Lambda_{\hg\hh}$-module of rank $1$. Then we have the following factorization of the image of the Beilinson--Flach class under the Perrin-Riou map $\Log_{\bm{\omega}_{\hg} \otimes \bm{\omega}_{\hh}}$: 
\[ \Omega_{f,\gamma} \cdot \Lp^g(f,\hg,\hh) \cdot \Lp^f(f, \hg, \hh) =  \lambda_{N_g}(\hg) \cdot \Log_{\omega_{\hg} \otimes \omega_{\hh}}(\BF(f, \hg, \hh)). \] 
\end{proposition}

\begin{proof}
By Corollary~\ref{corollary:comparison}, we have
\[
\lambda_{N_g}(\hg)\cdot\BF(f,\hg,\hh)=\Omega_{f,\gamma}\cdot \Lp^g(f,\hg,\hh) \cdot \kappa(f,\hg,\hh).
\]
The result then follows immediately by taking $\Log_{\bm{\omega}_{\hg} \otimes \bm{\omega}_{\hh}}$ and using Theorem~\ref{thm:reclawdiagonalcycles}. 
\end{proof}

\section{The $p$-exceptional case}\label{sec:expectional}

We keep the notations and assumptions in the previous section, and consider now the situation where there is a trivial zero for the $p$-adic $L$-function $\Lp^g(f,\hg,\hh)$, as well as for the $p$-adic $L$-functions $L_p(\hg,\hh)$ and $L_p(\hg,\hh\otimes\varepsilon_K)$, at the point $(y_0,z_0)$ corresponding to the modular forms $(g,h)$ arising from the vanishing of an Euler factor at $p$. In this situation, the statements in \S\ref{subsec:nonexceptional} become trivial, so we introduce \textit{improved} classes and \textit{improved} $p$-adic $L$-functions to remove the vanishing Euler factor. 

Hence, along this section we will work under the following assumption.

\begin{ass}
With the previous notations, it holds that \[ \alpha_g \beta_h = p^{\frac{l_0+m_0-3}{2}}. \] 
\end{ass}

Since $g$ and $h$ are $p$-ordinary and $p\nmid \mathrm{cond}(h)$, it follows from this assumption together with the Ramanujan--Petersson conjecture that $(l_0,m_0)=(2,1)$ and that $g$ has conductor $N_gp$. Note that this situation includes in particular the case where $g$ is the modular form attached to an elliptic curve over $\bb{Q}$ with multiplicative reduction at $p$, which is of great arithmetic interest.



We consider the following codimension-$1$ subvariety of $\cl{U}_{\hg}\times \cl{U}_{\hh}$:
\[
\cl{Y}=\lbrace (y,z)\in \cl{U}_{\hg}\times\cl{U}_{\hh}\;\vert\; \kappa_{\hg}(y)=\kappa_{\hh}(z)+1\rbrace.
\]
Let $\cl{C}$ be an irreducible component of $\cl{Y}$ containing the point $(y_0,z_0)$. We denote by $\cl{O}_{\cl{C}}$ the corresponding ring of global functions. Note that an arithmetic point $(y,z)\in \cl{C}\subset \cl{U}_\hg\times\cl{U}_\hh$ will have weights $(l,l-1)$ for some integer $l$. We denote by $\cl{C}^{\mathrm{cl}}$ the intersection of $\cl{C}$ with $\tilde{\cl{W}}_{\hg}^{\mathrm{cl}}\times \tilde{\cl{W}}_{\hh}^{\mathrm{cl}}$.

Note that, for any point $(y,z)\in \cl{C}^{\mathrm{cl}}$ with $p\nmid \mathrm{cond}(\hg_y)\cdot\mathrm{cond}(\hh_z)$, the Euler factor $\cl{E}(x_0,y,z)$ appearing in the interpolation formula for $\Lp^g(f,\hg,\hh)$ and the Euler factor $\cl{E}(y,z,s(y,z))$ appearing in the interpolation formula for $L_p(\hg,\hh,(\mathbf{l}+\mathbf{m}-1)/2)$ are given by
\[
\cl{E}(x_0,y,z)=\cl{E}(y,z,s(y,z))=\left(1-\chi_g(p)\frac{\alpha_{\hh_z}}{\alpha_{\hg_y}}\right)^2 \left(1-\chi_g(p)\frac{\beta_{\hh_z}}{\alpha_{\hg_y}}\right)^2.
\]
For varying values of $(y,z)\in \cl{C}^{\mathrm{cl}}$, the factors $1-\chi_g(p)\alpha_{\hh_z}\alpha_{\hg_y}^{-1}$ are interpolated by the $p$-adic $L$-function
\[
1-\chi_g(p)\frac{a_p(\hh)}{a_p(\hg)} \in \cl{O}_{\cl{C}},
\]
which vanishes at the point $(y_0,z_0)$.


\begin{remark}
The appearance of the double factor
\[
\left(1-\chi_g(p)\frac{a_p(\hh)}{a_p(\hg)}\right)^2
\]
in the interpolation formula for $\Lp^g(f,\hg,\hh)$ (resp. $L_p(\hg,\hh,(\mathbf{l}+\mathbf{m}-1)/2)$) at points in $\cl{C}^{\mathrm{cl}}$ comes from two different sources:
\begin{itemize}
    \item the numerator in the interpolation formula for the Perrin-Riou map $\cl{L}_{f\hg\hh}^{+-+}$ (resp. $\cl{L}_{\hg\hh}^{-+})$;
    \item the comparison between the corresponding $p$-stabilized and non-$p$-stabilized cohomology classes.
\end{itemize}
\end{remark}


\subsection{Improved $p$-adic $L$-functions}

In this subsection we introduce the improved $p$-adic $L$-functions that we will use.


\begin{theorem}\label{thm:improvedHRpadicLfunction}
    There exists a $p$-adic $L$-function
    \[
    \widehat{L}_p(\hg,\hh)\in I_{\hg}^{-1}\cl{O}_{\cl{C}}
    \]
    such that
    \[
    L_p(\hg,\hh,(\mathbf{l}+\mathbf{m}-1)/2)\vert_{\cl{C}}=\left(1-\chi_g(p)\frac{a_p(\hh)}{a_p(\hg)}\right)^2 \widehat{L}_p(\hg,\hh).
    \]
\end{theorem}




\begin{proof}
Let $N=\lcm(N_g,N_h)$. Let $E_1=E^{(1)}_{1/N}(\tau)$ be the weight-$1$ holomorphic Eisenstein series introduced in \cite[\S5.1]{LLZ0}. Then, as in the proof of \cite[Lemma~9.8]{BSV}, one shows that the $p$-adic $L$-function
\[
\widehat{L}_p(\hg,\hh)=\frac{\langle W_{N}(\hg), \mathrm{e}_{\ord}(E_1\cdot \hh)\rangle_N}{\langle \hg,\hg\rangle_N} \in I_{\hg}^{-1} \cl{O}_{\cl{C}}
\]
satisfies the required properties.

\end{proof}




\begin{theorem}\label{thm:improvedtripleproductpadicLfunction}
    There exists a $p$-adic $L$-function
    \[
    \widehat{\Lp}^g(f,\hg,\hh)\in I_{\hg}^{-1}\cl{O}_{\cl{C}}
    \]
    such that
    \[
    \Lp^g(f,\hg,\hh)\vert_{\cl{C}}=\left(1-\chi_g(p)\frac{a_p(\hh)}{a_p(\hg)}\right)^2 \widehat{\Lp}^g(f,\hg,\hh).
    \]
\end{theorem}
\begin{proof}
    This is proved as in \cite[Lemma~9.8]{BSV}, taking $(\hf,\hg,\hh)$ in \emph{loc.\,cit.} to be our $(\hg,\hh,\hf)$.
\end{proof}




\subsection{Improved Perrin-Riou maps}\label{subsec:improvedPerrinRiou}


In this subsection we introduce improved Perrin-Riou maps.  Given a $\Lambda_{\hg\hh}$-module $M$, we denote by $M\vert_{\cl{C}}$ the $\cl{O}_{\cl{C}}$-module $M\otimes_{\Lambda_{\hg\hh}}\cl{O}_{\cl{C}}$.

\begin{propo}\label{perrin-improved1}
There exists a homomorphism of $\mathcal{O}_{\cl{C}}$-modules
\[
\widehat{\mathcal L}_{\hg \hh}^{-+}: H^1(\mathbb{Q}_p,\mathbb{V}_{\hg}^- \hat{\otimes} \mathbb{V}_{\hh}^+ (2-\mathbf{t}_1)|_{\mathcal C}) \rightarrow I_\hg^{-1} \cl{O}_{\cl{C}},
\]
such that for every point $(y,z) \in \mathcal{C}^{\mathrm{cl}}$ of weights $(l,l-1)$ the specialization of $\widehat{\mathcal L}_{\hg \hh}^{-+}$ at $(y,z)$ is the homomorphism
\[
\widehat{\mathcal L}_{{\hg\hh}}^{-+}(y,z): H^1(\mathbb Q_p, V_{\hg_y}^- \otimes V_{\hh_z}^+(2-l)) \rightarrow \mathbb C_p
\]
given by
\[
\widehat{\mathcal L}_{{\hg \hh}}^{-+}(y,z) = \frac{1}{1-p^{1-l} \alpha_{\hg_y} \beta_{\hh_z}}  \times \langle \exp_{\BK}^*(\cdot), \omega_{f} \otimes\eta_{\hg_y}^\alpha \otimes \omega_{\hh_z} \rangle.
\]
\end{propo}
\begin{proof}
This follows from \cite[Lemma~9.4]{BSV}.
\end{proof}

\begin{propo}\label{perrin-improved2}
There exists a homomorphism of $\mathcal{O}_{\cl{C}}$-modules
\[
\widehat{\mathcal L}_{f \hg \hh}^{+-+}: H^1(\mathbb{Q}_p, V_f^+\otimes\mathbb{V}_{\hg}^- \hat{\otimes} \mathbb{V}_{\hh}^+ (2-\mathbf{t}_1)|_{\mathcal C}) \rightarrow I_\hg^{-1} \cl{O}_{\cl{C}},
\]
such that for every point $(y,z) \in \mathcal{C}^{\mathrm{cl}}$ of weights $(l,l-1)$ the specialization of $\widehat{\mathcal L}_{f\hg \hh}^{+-+}$ at $(y,z)$ is the homomorphism
\[
\widehat{\mathcal L}_{{f\hg\hh}}^{+-+}(y,z): H^1(\mathbb Q_p, V_f^+\otimes V_{\hg_y}^- \otimes V_{\hh_z}^+(2-l)) \rightarrow \mathbb C_p
\]
given by
\[
\widehat{\mathcal L}_{{f\hg \hh}}^{+-+}(y,z) = \frac{1}{1-p^{1-l} \alpha_{\hg_y} \beta_{\hh_z}}  \times \langle \exp_{\BK}^*(\cdot), \omega_{f} \otimes\eta_{\hg_y}^\alpha \otimes \omega_{\hh_z} \rangle.
\]
\end{propo}
\begin{proof}
This follows from \cite[Lemma~9.4]{BSV}.
\end{proof}

\begin{remark}
As before, given an element
\[
v\otimes z\in V_f^+\otimes H^1(\mathbb{Q}_p,\mathbb{V}_{\hg}^- \hat{\otimes} \mathbb{V}_{\hh}^+ (2-\mathbf{t}_1)|_{\mathcal C}) \cong H^1(\mathbb{Q}_p, V_f^+\otimes\mathbb{V}_{\hg}^- \hat{\otimes} \mathbb{V}_{\hh}^+ (2-\mathbf{t}_1)|_{\mathcal C}),
\]
we have that
\[
\widehat{\cl{L}}_{f\hg\hh}^{+-+}(v\otimes z)=\langle v,\omega_f\rangle\cdot \widehat{\cl{L}}_{\hg\hh}^{+-}(z).
\]
\end{remark}

\subsection{Improved cohomology classes}

We now introduce the improved version of the cohomology classes that have been used in this work. 

We can define Selmer conditions for the $\cl{O}_{\cl{C}}[G_{\bb{Q}}]$-modules $\bb{V}_{\hg\hh}^\dagger\vert_{\cl{C}}$ and $\bb{V}_{f\hg\hh}^\dagger\vert_{\cl{C}}$ analogous to the ones introduced for $\bb{V}_{\hg\hh}^\dagger$ and $\bb{V}_{f\hg\hh}^\dagger$ in the previous section and consider the corresponding Selmer groups. In particular, we can consider the Selmer groups $H^1_{\bal}(\bb{Q},\bb{V}_{\hg\hh}^\dagger\vert_{\cl{C}})$, $H^1_{\bal}(\bb{Q},\bb{V}_{f\hg\hh}^\dagger\vert_{\cl{C}})$ and $H^1_{\cl{G}\cup +}(\bb{Q},\bb{V}_{f\hg\hh}^\dagger\vert_{\cl{C}})$.

While in the case of diagonal cycles there is a construction of an improved cohomology class, this is not the case for Beilinson--Flach elements. As it occurs with the $p$-adic $L$-function, the factor $\left( 1 - \chi_g(p)\frac{a_p(\hh)}{a_p(\hg)} \right)$ appears in the interpolation property when considering its variation in families (see e.g. \cite[\S8]{KLZ}). Hence, the following conjecture can be seen as a standard expectation in the theory of exceptional zeros, and we expect to come back to it in forthcoming work.

\begin{conj}\label{conj:improvedbeilinsonclass}
There exists a cohomology class $\widehat{\kappa}_{\hg,\hh}\in H^1_{\bal}(\bb{Q},\bb{V}_{\hg\hh}^\dagger\vert_{\cl{C}})$ such that
\[
\kappa_{\hg,\hh}\vert_{\cl{C}} = \left( 1 - \chi_g(p)\frac{a_p(\hh)}{a_p(\hg)} \right) \widehat{\kappa}_{\hg,\hh}.
\]
\end{conj}

\begin{remark}
    Similarly, we also expect to have an analogous improved cohomology class $\widehat{\kappa}_{\hg,\hh\otimes\varepsilon_K} \in H^1_{\bal}(\bb{Q},\bb{V}_{\hg\hh}^\dagger\vert_{\cl{C}}\otimes\varepsilon_K)$.
\end{remark}

For the rest of this section, we work under the following assumption.

\begin{ass}
    Conjecture~\ref{conj:improvedbeilinsonclass} holds.
\end{ass}

Let $\mathrm{pr}^{-+}:H^1_{\bal}(\bb{Q}_p,\bb{V}_{\hg\hh}^\dagger\vert_{\cl{C}})\rightarrow H^1(\bb{Q}_p,\bb{V}_{\hg}^-\hat{\otimes}\bb{V}_\hh^+(2-\mathbf{t}_1)\vert_{\cl{C}})$ be the map induced by the natural projection $\mathscr{F}^2\bb{V}_{\hg\hh}^\dagger\vert_{\cl{C}}\rightarrow \bb{V}_\hg^-\hat{\otimes}\bb{V}_\hh^+(2-\mathbf{t}_1)\vert_{\cl{C}}$ and let
\[
\widehat{\mathcal{L}}_{\hg\hh}^{\hg}:H^1_{\bal}(\bb{Q},\bb{V}_{\hg\hh}^\dagger\vert_{\cl{C}})\longrightarrow I_{\hg}^{-1}\cl{O}_{\cl{C}}
\]
be the map defined by $\widehat{\mathcal{L}}_{\hg\hh}^\hg=\widehat{\mathcal{L}}_{\hg\hh}^{-+}\circ \mathrm{pr}^{-+}\circ \res_p$.

\begin{propo}\label{prop:improvedreclaw1}
    The class $\widehat{\kappa}_{\hg,\hh}$ satisfies
    \[
    \widehat{\cl{L}}_{\hg\hh}^{\hg}(\widehat{\kappa}_{\hg,\hh})=\widehat{L}_p(\hg,\hh).
    \]
\end{propo}
\begin{proof}
   By Proposition~\ref{perrin-improved1}, Theorem~\ref{reclaw} and Theorem~\ref{thm:improvedHRpadicLfunction}, we have
    \begin{align*}
        \left( 1 - \chi_g(p)\frac{a_p(\hh)}{a_p(\hg)} \right)^2  \widehat{\cl{L}}_{\hg\hh}^\hg (\widehat{\kappa}_{\hg,\hh})&=\cl{L}_{\hg\hh}^{\hg}(\kappa_{\hg,\hh})\vert_{\cl{C}}= L_p(\hg,\hh)\vert_{\cl{C}} \\
        &=\left(1-\chi_g(p)\frac{a_p(\hh)}{a_p(\hg)}\right)^2 \widehat{L}_p(\hg,\hh).
    \end{align*}
    By the Ramanujan--Petersson conjecture, the factor $1-\chi_g(p)\frac{a_p(\hh)}{a_p(\hg)}$ does not vanish at any point $(y,z)\in\cl{C}^{\mathrm{cl}}$ with $p\nmid \mathrm{cond}(\hg_y)\cdot \mathrm{cond}(\hh_z)$, so in particular it is non-zero in $\cl{O}_{\cl{C}}$. The result now follows from the above equality. 
\end{proof}

\begin{remark}
    Similarly $\widehat{\cl{L}}_{\hg\hh}^{\hg}(\widehat{\kappa}_{\hg,\hh\otimes\varepsilon_K})=\widehat{L}_p(\hg,\hh\otimes\varepsilon_K)$.
\end{remark}


\begin{propo}\label{prop:improvedtripleproductclass}
There exists a cohomology class $\widehat{\kappa}(f,\hg,\hh)\in H^1_{\bal}(\bb{Q},\bb{V}_{f\hg\hh}^\dagger\vert_{\cl{C}})$ such that
\[
\kappa(f,\hg,\hh)\vert_{\cl{C}} = \left( 1 - \chi_g(p)\frac{a_p(\hh)}{a_p(\hg)} \right) \widehat{\kappa}(f,\hg,\hh).
\]
\end{propo}
\begin{proof}
This is proved in \cite[\S9.3]{BSV}.
\end{proof}

Let $\mathrm{pr}^{+-+}:H^1_{\bal}(\bb{Q}_p,\bb{V}_{f\hg\hh}^\dagger\vert_{\cl{C}})\rightarrow H^1(\bb{Q}_p,V_f^+\otimes\bb{V}_{\hg}^-\hat{\otimes}\bb{V}_\hh^+(2-\mathbf{t}_1)\vert_{\cl{C}})$ be the map induced by the natural projection $\mathscr{F}^2\bb{V}_{f\hg\hh}^\dagger\vert_{\cl{C}}\rightarrow V_f^+\otimes\bb{V}_\hg^-\hat{\otimes}\bb{V}_\hh^+(2-\mathbf{t}_1)\vert_{\cl{C}}$ and let
\[
\widehat{\mathcal{L}}_{f\hg\hh}^{\hg}:H^1_{\bal}(\bb{Q},\bb{V}_{f\hg\hh}^\dagger\vert_{\cl{C}})\longrightarrow I_{\hg}^{-1}\cl{O}_{\cl{C}}
\]
be the map defined by $\widehat{\mathcal{L}}_{f\hg\hh}^\hg=\widehat{\mathcal{L}}_{f\hg\hh}^{+-+}\circ \mathrm{pr}^{+-+}\circ \res_p$.

\begin{propo}
    The class $\widehat{\kappa}(f,\hg,\hh)$ satisfies
    \[
    \widehat{\cl{L}}_{f\hg\hh}^{\hg}(\widehat{\kappa}(f,\hg,\hh))=\widehat{\Lp}^g(f,\hg,\hh).
    \]
\end{propo}
\begin{proof}
    By Proposition~\ref{perrin-improved2}, Proposition~\ref{prop:improvedtripleproductclass}, Theorem~\ref{thm:reclawdiagonalcycles} and Theorem~\ref{thm:improvedtripleproductpadicLfunction}, we have
    \begin{align*}
        \left( 1 - \chi_g(p)\frac{a_p(\hh)}{a_p(\hg)} \right)^2  \widehat{\cl{L}}_{f\hg\hh}^\hg (\widehat{\kappa}(f,\hg,\hh))&=\cl{L}_{f\hg\hh}^{\hg}(\kappa(f,\hg,\hh))\vert_{\cl{C}}= \Lp^g(f,\hg,\hh)\vert_{\cl{C}} \\
        &=\left(1-\chi_g(p)\frac{a_p(\hh)}{a_p(\hg)}\right)^2 \widehat{\Lp}(f,\hg,\hh).
    \end{align*}
    By the Ramanujan--Petersson conjecture, the factor $1-\chi_g(p)\frac{a_p(\hh)}{a_p(\hg)}$ does not vanish at any point $(y,z)\in\cl{C}^{\mathrm{cl}}$ with $p\nmid \mathrm{cond}(\hg_y)\cdot \mathrm{cond}(\hh_z)$, so in particular it is non-zero in $\cl{O}_{\cl{C}}$. The result now follows from the above equality.
\end{proof}

\subsection{The main theorem}

In this subsection, we adapt the main result of \S\ref{subsec:nonexceptional} to the current exceptional case.

\begin{propo}\label{prop:improvedfactorization}
    We have the following equality of $p$-adic $L$-functions:
    \[
    \widehat{\Lp}^g(f,\hg,\hh)^2=\widehat{L}_p(\hg,\hh)\widehat{L}_p(\hg,\hh\otimes\varepsilon_K).
    \]
\end{propo}
\begin{proof}
    By Proposition~\ref{prop:factorization}, Theorem~\ref{thm:improvedHRpadicLfunction} and Theorem~\ref{thm:improvedtripleproductpadicLfunction}, we have the equality
    \[
    \left(1-\chi_g(p)\frac{a_p(\hh)}{a_p(\hg)}\right)^4 \widehat{\Lp}^g(f,\hg,\hh)^2=\left(1-\chi_g(p)\frac{a_p(\hh)}{a_p(\hg)}\right)^4\widehat{L}_p(\hg,\hh)\widehat{L}_p(\hg,\hh\otimes\varepsilon_K)
    \]
    in $I_{\hg}^{-2}\cl{O}_{\cl{C}}$. By the Ramanujan--Petersson conjecture, the factor $1-\chi_g(p)\frac{a_p(\hh)}{a_p(\hg)}$ does not vanish at any point $(y,z)\in\cl{C}^{\mathrm{cl}}$ with $p\nmid \mathrm{cond}(\hg_y)\cdot \mathrm{cond}(\hh_z)$, so in particular it is non-zero in $\cl{O}_{\cl{C}}$. The result now follows from the above equality.
\end{proof}

We now introduce the improved Beilinson--Flach class that we will use. For the remaining of the section, we assume that Conjecture~\ref{conj:improvedbeilinsonclass} is true.

\begin{defi}
The \emph{improved weighted Beilinson--Flach class} associated with the pair $(\hg,\hh)$ is the element
\[
\widehat{\mathrm{BF}}(f,\hg,\hh)=v_{f,1}\otimes \widehat{L}_p(\hg,\hh\otimes\varepsilon_K)\widehat{\kappa}_{\hg,\hh}+v_{f,\varepsilon_K}\otimes\widehat{L}_p(\hg,\hh)\widehat{\kappa}_{\hg,\hh\otimes \varepsilon_K}
\]
in $H^1(\bb{Q},\bb{V}_{f\hg\hh}^\dagger\vert_{\cl{C}})$.   
\end{defi}

Let $\res_p(\widehat{\kappa}(f,\hg,\hh))^{+-+}$ and $\res_p(\widehat{\BF}(f,\hg,\hh))^{+-+}$ be the images of $\res_p(\widehat{\kappa}(f,\hg,\hh))$ and $\res_p(\widehat{\BF}(f,\hg,\hh))$, respectively, in $H^1(\bb{Q}_p,V_f^+\otimes\bb{V}_\hg^-\hat{\otimes}\bb{V}_\hh^+(2-\mathbf{t}_1)\vert_{\cl{C}})$.

\begin{theorem}\label{thm:bf-vs-cycles-fam}
The equality
\[
\Omega_{f,\gamma} \cdot \widehat{\cl{L}}_{f\hg\hh}^{+-+}(\res_p(\widehat{\kappa}(f,\hg,\hh))^{+-+})^2 =\lambda_{N_g}(\hg) \cdot \widehat{\mathcal{L}}_{f\hg\hh}^{+-+}(\res_p(\widehat{\BF}(f, \hg, \hh))^{+-+})
\]
holds. 
\end{theorem}
\begin{proof}
The result follows as in Theorem~\ref{thm:main}, working in this case with the improved Perrin-Riou maps.
\end{proof}

\begin{proposition}
    Assume that $H^1(\bb{Q},\bb{V}_{f\hg\hh}^\dagger\vert_{\cl{C}})$ is a torsion-free $\cl{O}_{\cl{C}}$-module and $H^1_{\cl{G}\cup +}(\bb{Q},\bb{V}_{f\hg\hh}^\dagger)$ is a torsion-free $\Lambda_{\hg\hh}$-module of rank $1$. Then $\widehat{\mathrm{BF}}(f, \hg, \hh)$ belongs to $H^1_{\bal}(\bb{Q},\bb{V}_{f\hg\hh}^\dagger\vert_{\cl{C}})$ and
    \[
    \Omega_{f,\gamma} \cdot \widehat{\Lp}^g(f,\hg,\hh) \cdot \widehat{\kappa}(f,\hg,\hh) = \lambda_{N_g}(\hg)\cdot\widehat{\mathrm{BF}}(f, \hg, \hh).
    \]
\end{proposition}
\begin{proof}
    By Corollary~\ref{corollary:comparison}, we have
    \[
    \left(1-\chi_g(p)\frac{a_p(\hh)}{a_p(\hg)}\right)^3\cdot\Omega_{f,\gamma} \cdot \widehat{\Lp}^g(f,\hg,\hh) \cdot \widehat{\kappa}(f,\hg,\hh) = \left(1-\chi_g(p)\frac{a_p(\hh)}{a_p(\hg)}\right)^3\cdot \lambda_{N_g}(\hg)\cdot\widehat{\mathrm{BF}}(f, \hg, \hh).
    \]
    Since $H^1(\bb{Q},\bb{V}_{f\hg\hh}^\dagger\vert_{\cl{C}})$ is a torsion-free $\cl{O}_{\cl{C}}$-module and $1-\chi_g(p)\frac{a_p(\hh)}{a_p(\hg)}$ is non-zero in $\cl{O}_{\cl{C}}$, it follows that
    \[
    \Omega_{f,\gamma} \cdot \widehat{\Lp}^g(f,\hg,\hh) \cdot \widehat{\kappa}(f,\hg,\hh) = \lambda_{N_g}(\hg)\cdot\widehat{\mathrm{BF}}(f, \hg, \hh).
    \]
    In particular, since $\widehat{\kappa}(f,\hg,\hh)$ belongs to $H^1_{\bal}(\bb{Q},\bb{V}_{f\hg\hh}^\dagger\vert_{\cl{C}})$, the same is true for $\widehat{\BF}(f,\hg,\hh)$.
\end{proof}

\begin{remark}
    As in Remark~\ref{rk:torsionfreeness}, we can ensure that both $H^1_{\cl{G}\cup +}(\bb{Q}, \bb{V}_{f\hg\hh}^\dagger)$ and $H^1(\bb{Q},\bb{V}_{f\hg\hh}^\dagger\vert_{\cl{C}})$ are torsion-free by imposing the condition that $H^0(\bb{Q},\overline{\rho}^\dagger)=0$, where $\overline{\rho}^\dagger$ is the residual representation attached to $\bb{V}_{f\hg \hh}^\dagger$.
\end{remark}

Let $\widehat{\kappa}(f,g_\alpha,h_\alpha)\in H^1_{\bal}(\bb{Q},V_{fgh}^\dagger)$ and $\widehat{\BF}(f,g_\alpha,h_\alpha)\in H^1_{\cl{G}\cup +}(\bb{Q},V_{fgh}^\dagger)$ be the specializations of the classes $\widehat{\kappa}(f,\hg,\hh)\in H^1_{\bal}(\bb{Q},\bb{V}_{f\hg\hh}^\dagger\vert_{\cl{C}})$ and $\widehat{\BF}(f,\hg,\hh)\in H^1(\bb{Q},\bb{V}_{f\hg\hh}^\dagger\vert_{\cl{C}})$, respectively, at the point $(y_0,z_0)$.

\begin{corollary}
    Assume that $H^1(\bb{Q},\bb{V}_{f\hg\hh}^\dagger\vert_{\cl{C}})$ is a torsion-free $\cl{O}_{\cl{C}}$-module and $H^1_{\cl{G}\cup +}(\bb{Q},\bb{V}_{f\hg\hh}^\dagger)$ is a torsion-free $\Lambda_{\hg\hh}$-module of rank $1$. Then $\widehat{\BF}(f,g_\alpha,h_\alpha)$ belongs to $H^1_{\bal}(\bb{Q}, V_{fgh}^\dagger)$ and
    \[
    \lambda_{N_g}(g)\cdot\widehat{\BF}(f,g_\alpha,h_\alpha)=\Omega_{f,\gamma}\cdot \widehat{\Lp}^g(f,g_\alpha,h_\alpha) \cdot \widehat{\kappa}(f,g_\alpha,h_\alpha).
    \]
\end{corollary}
\begin{proof}
    This is an immediate consequence of the previous proposition.
\end{proof}

\bibliographystyle{amsalpha}
\bibliography{AOR-refs}

\providecommand{\bysame}{\leavevmode\hbox to3em{\hrulefill}\thinspace}
\providecommand{\MR}{\relax\ifhmode\unskip\space\fi MR }
\providecommand{\MRhref}[2]{%
  \href{http://www.ams.org/mathscinet-getitem?mr=#1}{#2}
}
\providecommand{\href}[2]{#2}
\begin{thebibliography}{BSTW24}

\bibitem[ACR23a]{ACR2}
Ra\'ul Alonso, Francesc Castella, and \'Oscar Rivero, \emph{An anticyclotomic euler system for for adjoint modular galois representations}, Ann. Inst. Fourier (Grenoble), to appear (2023), available at \href{https://arxiv.org/abs/2204.07658}{arXiv:2204.07658}.

\bibitem[ACR23b]{ACR}
\bysame, \emph{The diagonal cycle euler system for {$\GL_2\times \GL_2$}}, J. Inst. Math. Jussieu, to appear (2023), available at \href{https://arxiv.org/abs/2106.05322}{arXiv:2106.05322}.

\bibitem[AI21]{AI}
Fabrizio Andreatta and Adrian Iovita, \emph{Triple product {$p$}-adic {$L$}-functions associated to finite slope {$p$}-adic families of modular forms}, Duke Math. J. \textbf{170} (2021), no.~9, 1989--2083. \MR{4278669}

\bibitem[BC23]{BC}
K\^azim B\"uy\"yukboduk and Daniele Casazza, \emph{On the {A}rtin formalism for triple product {$p$}-adic {$L$}-functions: Super-factorization}, preprint, \href{https://arxiv.org/abs/2301.08383}{arXiv:2301.08383}.

\bibitem[BD16]{bellaiche-dimitrov}
Jo\"el Bella\"iche and Mladen Dimitrov, \emph{On the eigencurve at classical weight 1 points}, Duke Math. J. \textbf{165} (2016), no.~2, 245--266. \MR{3457673}

\bibitem[BDP22]{BeDiPo}
Adel Betina, Mladen Dimitrov, and Alice Pozzi, \emph{On the failure of {G}orensteinness at weight 1 {E}isenstein points of the eigencurve}, Amer. J. Math. \textbf{144} (2022), no.~1, 227--265. \MR{4367418}

\bibitem[BDV22]{BDV}
Massimo Bertolini, Henri Darmon, and Rodolfo Venerucci, \emph{Heegner points and {B}eilinson-{K}ato elements: a conjecture of {P}errin-{R}iou}, Adv. Math. \textbf{398} (2022), Paper No. 108172, 50. \MR{4379203}

\bibitem[Bel12]{bellaiche11a}
Jo{\"e}l Bella{\"i}che, \emph{Critical {$p$}-adic {$L$}-functions}, Invent. Math. \textbf{189} (2012), no.~1, 1--60. \MR{2929082}

\bibitem[BS23]{BS}
K\^azim B\"uy\"yukboduk and Ryotaro Sakamoto, \emph{On the {A}rtin formalism for triple product {$p$}-adic {$L$}-functions}, preprint, \href{https://arxiv.org/abs/2501.06541}{arXiv:2501.06541}.

\bibitem[BSTW24]{BSTW}
Ashay Burungale, Christopher Skinner, Ye~Tian, and Xin Wan, \emph{Zeta elements for elliptic curves and applications}, preprint, {\tt arXiv:2409.01350} (2024).

\bibitem[BSV22a]{BSV2}
Massimo Bertolini, Marco~Adamo Seveso, and Rodolfo Venerucci, \emph{Balanced diagonal classes and rational points on elliptic curves}, no. 434, 2022, Heegner points, Stark-Heegner points, and diagonal classes, pp.~175--201. \MR{4489474}

\bibitem[BSV22b]{BSV}
\bysame, \emph{Reciprocity laws for balanced diagonal classes}, no. 434, 2022, Heegner points, Stark-Heegner points, and diagonal classes, pp.~77--174. \MR{4489473}

\bibitem[Buy16]{buyukboduk16}
K\^{a}z\i{m} Buy\"ukboduk, \emph{On {N}ekov\'a\v r's heights, exceptional zeros and a conjecture of {M}azur-{T}ate-{T}eitelbaum}, Int. Math. Res. Not. IMRN (2016), no.~7, 2197--2237. \MR{3509952}

\bibitem[CD23]{castella-do}
Francesc Castella and Kim~Tuan Do, \emph{Diagonal cycles and anticyclotomic iwasawa theory of modular forms}, preprint, \href{https://arxiv.org/abs/2303.06751}{arXiv:2303.06751}.

\bibitem[DR14]{DR1}
Henri Darmon and Victor Rotger, \emph{Diagonal cycles and {E}uler systems {I}: {A} {$p$}-adic {G}ross-{Z}agier formula}, Ann. Sci. \'Ec. Norm. Sup\'er. (4) \textbf{47} (2014), no.~4, 779--832.

\bibitem[DR17]{DR2}
\bysame, \emph{Diagonal cycles and {E}uler systems {II}: {T}he {B}irch and {S}winnerton-{D}yer conjecture for {H}asse-{W}eil-{A}rtin {$L$}-functions}, J. Amer. Math. Soc. \textbf{30} (2017), no.~3, 601--672.

\bibitem[DR22]{DR3}
\bysame, \emph{{$p$}-adic families of diagonal cycles}, no. 434, 2022, Heegner points, Stark-Heegner points, and diagonal classes, pp.~29--75. \MR{4489472}

\bibitem[Hid85]{hida1}
Haruzo Hida, \emph{A {$p$}-adic measure attached to the zeta functions associated with two elliptic modular forms. {I}}, Invent. Math. \textbf{79} (1985), no.~1, 159--195. \MR{774534}

\bibitem[Hid88]{hida2}
\bysame, \emph{A {$p$}-adic measure attached to the zeta functions associated with two elliptic modular forms. {II}}, Ann. Inst. Fourier (Grenoble) \textbf{38} (1988), no.~3, 1--83. \MR{976685}

\bibitem[Hsi21]{Hs}
Ming-Lun Hsieh, \emph{Hida families and $p$-adic triple product {$L$}-functions}, American Journal of Mathematics \textbf{143} (2021), no.~2, 411--532.

\bibitem[KLZ17]{KLZ}
Guido Kings, David Loeffler, and Sarah~Livia Zerbes, \emph{Rankin-{E}isenstein classes and explicit reciprocity laws}, Camb. J. Math. \textbf{5} (2017), no.~1, 1--122.

\bibitem[KLZ20]{KLZ0}
\bysame, \emph{Rankin-{E}isenstein classes for modular forms}, Amer. J. Math. \textbf{142} (2020), no.~1, 79--138.

\bibitem[LLZ14]{LLZ0}
Antonio Lei, David Loeffler, and Sarah~Livia Zerbes, \emph{Euler systems for {R}ankin-{S}elberg convolutions of modular forms}, Ann. of Math. (2) \textbf{180} (2014), no.~2, 653--771.

\bibitem[LR24]{LR1}
David Loeffler and \'Oscar Rivero, \emph{Eisenstein degeneration of {E}uler systems}, J. Reine Angew. Math. \textbf{814} (2024), 241--282. \MR{4793345}

\bibitem[LZ16]{LZ-coleman}
David Loeffler and Sarah~Livia Zerbes, \emph{Rankin-{E}isenstein classes in {C}oleman families}, Res. Math. Sci. \textbf{3} (2016), Paper No. 29, 53.

\bibitem[Poz19]{pozzi-thesis}
Alice Pozzi, \emph{The {E}igencurve at {W}eight {O}ne {E}isenstein {P}oints}, Ph.D. thesis, McGill University, Montr\'eal, 2019. \MR{4197785}

\bibitem[PR25]{PR25}
Javier Polo and \'Oscar Rivero, \emph{Eisenstein degeneration of beilinson--kato classes and circular units}, preprint, 2025.

\bibitem[Urb14]{Urban-nearly-overconvergent}
Eric Urban, \emph{Nearly overconvergent modular forms}, Iwasawa Theory 2012: State of the Art and Recent Advances (Berlin), Contrib. Math. Comput. Sci., vol.~7, Springer, 2014, pp.~401--441.

\bibitem[Ven16]{Ven}
Rodolfo Venerucci, \emph{Exceptional zero formulae and a conjecture of {P}errin-{R}iou}, Invent. Math. \textbf{203} (2016), no.~3, 923--972. \MR{3461369}

\end{thebibliography}

\end{document}